\documentclass[reqno]{lms}
\usepackage{diagrams,amssymb,amsmath}

\newcommand{\mynote}[1]{\noindent\textbf{[#1]}}

\diagramstyle[height=2em,width=2em,midshaft,labelstyle=\scriptstyle]
\newarrow{To}----{->}
\newarrow{Epi}----{->>}
\newarrow{Mono}>---{->}
\newarrow{Iso}>---{->>}
\newarrow{Mapsto}|---{->}
\newarrow{Igual}=====
\newarrow{Dashto}{}{dash}{}{dash}{->}

\newcommand{\Z}{{\mathbb Z}}

\newcommand{\F}{{\mathbb F}}

\newcommand{\thetaS}{f}   
\newcommand{\thetaSb}{f}  

\newcommand{\pcom}{{}_{p}^{\wedge}}

\newcommand{\zploc}{\Z_{(p)}}

\DeclareMathAlphabet\EuR{U}{eur}{m}{n}
\SetMathAlphabet\EuR{bold}{U}{eur}{b}{n}

\newcommand{\Inj}{\operatorname{Inj}\nolimits}

\newcommand{\defeq}{\overset{\text{\textup{def}}}{=}}
\newcommand{\gen}[1]{{\langle}#1{\rangle}}

\renewcommand{\:}{\colon}

\newcommand{\calb}{\mathcal{B}}

\newcommand{\calf}{\mathcal{F}}

\newcommand{\call}{\mathcal{L}}

\newcommand{\calp}{\mathcal{P}}

\newcommand{\cals}{\mathcal{S}}
\newcommand{\calz}{\mathcal{Z}}

\newcommand{\orb}{\mathcal{O}}

\newcommand{\nv}[1]{|#1|} 

\newcommand{\widebar}[1]{\overset{\mskip3mu\hrulefill\mskip3mu}{#1}
                \vphantom{#1}}
\newcommand{\sminus}{\smallsetminus}

\newcommand{\nsg}{\vartriangleleft}
\newcommand{\Id}{\textup{Id}}
\newcommand{\incl}{\operatorname{incl}\nolimits}

\newcommand{\Inn}{\textup{Inn}}

\let\oldcirc=\circ
\renewcommand{\circ}{\mathchoice
    {\mathbin{\scriptstyle\oldcirc}}{\mathbin{\scriptstyle\oldcirc}}
    {\mathbin{\scriptscriptstyle\oldcirc}}
    {\mathbin{\scriptscriptstyle\oldcirc}}}

\newcommand{\hclim}[1]{\setbox1=\hbox{\rm hocolim}
    \setbox2=\hbox to \wd1{\rightarrowfill} \ht2=0pt \dp2=-1pt
    \mathop{\vtop{\baselineskip=5pt\box1\box2}}
    _{#1}}

\newcommand{\map}{\operatorname{Map}\nolimits}

\newcommand{\Hom}{\operatorname{Hom}\nolimits}

\newcommand{\Iso}{\operatorname{Iso}\nolimits}
\newcommand{\Aut}{\operatorname{Aut}\nolimits}
\newcommand{\Out}{\operatorname{Out}\nolimits}

\newcommand{\Ob}{\operatorname{Ob}\nolimits}
\newcommand{\Mor}{\operatorname{Mor}\nolimits}
\newcommand{\Ker}{\operatorname{Ker}\nolimits}

\renewcommand{\Im}{\operatorname{Im}\nolimits}

\newcommand{\longleft}[1]{\;{\leftarrow%
\count255=0 \loop \mathrel{\mkern-6mu}%
    \relbar\advance\count255 by1\ifnum\count255<#1\repeat}\;}
\newcommand{\longright}[1]{\;{\count255=0 \loop \relbar\mathrel{\mkern-6mu}%
    \advance\count255 by1\ifnum\count255<#1\repeat\rightarrow}\;}
\newcommand{\Right}[2]{\overset{#2}{\longright#1}}
\newcommand{\RIGHT}[3]{\mathrel{\mathop{\kern0pt\longright#1}
        \limits^{#2}_{#3}}}

\newcommand{\LEFT}[3]{\mathrel{\mathop{\kern0pt\longleft#1}\limits^{#2}_{#3}}
}
\newcommand{\dRIGHT}[3]{\mathrel{%
   \mathop{\vcenter{\baselineskip=0pt\hbox{$\kern0pt\longright#1$}%
   \hbox{$\kern0pt\longright#1$}}}\limits^{#2}_{#3}}}
\newcommand{\LRIGHT}[3]{\mathrel{%
   \mathop{\vcenter{\baselineskip=0pt\hbox{$\kern0pt\longleft#1$}%
   \hbox{$\kern0pt\longright#1$}}}\limits^{#2}_{#3}}}
\newcommand{\RLEFT}[3]{\mathrel{%
   \mathop{\vcenter{\baselineskip=0pt\hbox{$\kern0pt\longright#1$}%
   \hbox{$\kern0pt\longleft#1$}}}\limits^{#2}_{#3}}}
\newcommand{\onto}[1]{\;{\count255=0 \loop \relbar\joinrel
    \advance\count255 by1
    \ifnum\count255<#1 \repeat \twoheadrightarrow}\;}

\newtheorem{Thm}{Theorem}[section]
\newtheorem{Prop}[Thm]{Proposition}
\newtheorem{Cor}[Thm]{Corollary}
\newtheorem{Lem}[Thm]{Lemma}

\newtheorem{Defi}[Thm]{Definition}
\newtheorem{Th}{Theorem}

\newtheorem{Prp}[Th]{Proposition}

\newcommand{\longline}{\bigskip\centerline{\hbox to 5cm{\hrulefill}}\bigskip}

\newcommand{\frc}[1]{#1^c}    

\newcommand{\cj}[3][]{\def\test{#1}\def\tst{x}\ifx\test\tst{#2^{-1}#3#2}
	\else{#2#3#2^{-1}}\fi}
\newcommand{\cjup}[3][]{\def\test{#1}\def\tst{x}\def\tset{-}
	\ifx\test\tst{{}^{#2^{-1}}\!#3} 
	\else{\ifx\test\tset{{}^{#2}#3} \else{{}^{#2}\!#3} \fi} \fi}

\renewcommand{\b}[1]{\check{#1}}

\newcommand{\higherlim}[2]{\displaystyle\setbox1=\hbox{\rm lim}
	\setbox2=\hbox to \wd1{\leftarrowfill} \ht2=0pt \dp2=-1pt
	\setbox3=\hbox{$\scriptstyle{#1}$}
	\def\test{#1}\ifx\test\empty
	\mathop{\mathop{\vtop{\baselineskip=5pt\box1\box2}}}\nolimits^{#2}
	\else
	\ifdim\wd1<\wd3
	\mathop{\hphantom{^{#2}}\vtop{\baselineskip=5pt\box1\box2}^{#2}}_{#1}
	\else
	\mathop{\mathop{\vtop{\baselineskip=5pt\box1\box2}}_{#1}}%
	\nolimits^{#2}
	\fi\fi}

\newcommand{\calh}{\mathcal{H}}

\newcommand{\qc}{^{q}}
\newcommand{\calfq}{\calf\qc}
\newcommand{\callq}{\call\qc}

\newcommand{\SFL}[1][]{(S#1,\calf#1,\call#1)}

\newcommand{\calt}{\mathcal{T}}

\renewcommand{\labelenumi}{\textup{(\alph{enumi})}}%

\newcommand{\homf}{\Hom_{\calf}}

\newcommand{\autf}{\Aut_{\calf}}
\newcommand{\outf}{\Out_{\calf}}
\newcommand{\isof}{\Iso_{\calf}}

\newcommand{\sylp}[1]{\textup{Syl}_p(#1)}
\let\out=\Out

\renewenvironment{enumerate}{\begin{list}%
{\labelenumi}
{\usecounter{enumi}%
\setlength{\itemindent}{0pt}%
\settowidth{\labelwidth}{\labelenumi}%
\addtolength{\labelwidth}{\labelsep}%
\setlength{\leftmargin}{\labelsep}%
\addtolength{\leftmargin}{\labelwidth}%
\setlength{\listparindent}{0pt}%
\setlength{\itemsep}{6pt}%
\setlength{\parsep}{0pt}%
\setlength{\topsep}{6pt}%
}}{\end{list}}

\newenvironment{enumtwo}{\begin{list}%
{\labelenumii}
{\usecounter{enumii}%
\setlength{\itemindent}{0pt}%
\settowidth{\labelwidth}{\labelenumii}%
\addtolength{\labelwidth}{\labelsep}%
\setlength{\leftmargin}{\labelsep}%
\addtolength{\leftmargin}{\labelwidth}%
\setlength{\listparindent}{0pt}%
\setlength{\itemsep}{6pt}%
\setlength{\parsep}{0pt}%
\setlength{\topsep}{6pt}%
}}{\end{list}}

\newcommand{\hilim}[4][]{\def\test{#1}\def\tst{-}
\ifx\test\tst{\setbox1=\hbox{\rm lim}
\setbox2=\hbox to \wd1{\leftarrowfill} \ht2=0pt \dp2=-1pt
\mathop{\vtop{\baselineskip=5pt\box1\box2}}\nolimits^{#3}(#4)}
\else{\higherlim{#2}{#3}(#4)}\fi}
\let\2=\partial
\def\Syl{\textup{Syl}}

\newcommand{\dirlim}[1]{\setbox1=\hbox{\rm colim}
	\setbox2=\hbox to \wd1{\rightarrowfill} \ht2=0pt \dp2=-1pt
	\mathop{\vtop{\baselineskip=5pt\box1\box2}}
	_{#1}}

\renewcommand{\onto}[1]{\;{\count255=0 \loop \relbar\mathrel{\mkern-6mu}%
	\advance\count255 by1 \ifnum\count255<#1\repeat\twoheadrightarrow}\;}

\def\beq#1\eeq{\begin{equation*}#1\end{equation*}}

\title{Subgroup families controlling $p$-local finite groups}%

\author{C. Broto, N. Castellana, J. Grodal, R. Levi, and B. Oliver}

\classno{Primary 55R35. Secondary 55R40, 20D20}

\extraline{
C. Broto and N. Castellana were partially supported by MCYT grant 
BFM2001--2035\\
J. Grodal was partially supported by NSF grant DMS-0104318\\
R. Levi was partially supported by EPSRC grant GR/M7831\\
B. Oliver was partially supported by UMR 7539 of the CNRS}

\begin{document}

\maketitle

\begin{abstract}  
A $p$-local finite group consists of a finite $p$-group $S$, together with 
a pair of categories which encode ``conjugacy'' relations among subgroups 
of $S$, and which are modelled on the fusion in a Sylow $p$-subgroup of a 
finite group.  It contains enough information to define a 
classifying space which has many of the same properties as $p$-completed  
classifying spaces of finite groups.  In this paper, we examine which 
subgroups control this structure.  More precisely, we prove that the 
question of whether an abstract fusion system $\calf$ over a finite 
$p$-group $S$ is saturated can be determined by just looking at smaller 
classes of subgroups of $S$.  We also prove that the homotopy type of the 
classifying space of a given $p$-local finite group is independent of the 
family of subgroups used to define it, in the sense that it remains 
unchanged when that family ranges from the set of $\calf$-centric 
$\calf$-radical subgroups (at a minimum) to the set of 
$\calf$-quasicentric subgroups (at a maximum).  Finally, we look at 
constrained fusion systems, analogous to $p$-constrained finite groups, 
and prove that they in fact all arise from groups.
\end{abstract}


A $p$-local finite group consists of a finite $p$-group S, together with
a pair of categories $(\calf,\call)$, of which $\calf$ is modeled on the 
conjugacy (or fusion) in a Sylow subgroup of a finite group. The category 
$\call$ is essentially an extension of $\calf$ and contains just enough
extra information so that its $p$-completed nerve has many of the same 
properties as $p$-completed classifying spaces of finite groups. The main 
purpose of this paper is to study when the set of objects of $\calf$ or 
$\call$ can be changed without changing the conjugacy encoded by $\calf$ 
or the homotopy type of the nerve of $\call$. The tools introduced 
simplify the construction and manipulation of $p$-local finite groups in 
many cases.

We first recall the fusion and linking categories associated to a finite 
group.  Fix a prime $p$, a finite group $G$, and a Sylow $p$-subgroup $S$ 
of $G$. A \emph{$p$-fusion category} for $G$ is a category 
$\calf=\calf_S^\calh(G)$, whose object set is a set $\calh$ of subgroups 
of $S$, and whose morphisms are the homomorphisms between subgroups in 
$\calh$ induced by conjugation in $G$. The associated linking category 
$\call=\call_S^\calh(G)$ has the same objects, and morphisms from $P$ to 
$Q$ are given by the formula
	$$ \Mor_\call(P,Q) = \{x\in{}G\,|\,xPx^{-1}\le Q\}/O^p(C_G(P)). $$
Here, $O^p(-)$ is the subgroup generated by elements of order prime to 
$p$.  There is a canonical quotient functor $\call_S^\calh(G)\Right2{} 
\calf_S^\calh(G)$ which sends the class of $x$ to conjugation by $x$.  It 
was shown in \cite{BLO1} that the homotopy theory of the nerve 
$|\call^{\calh}_S(G)|$ (for the right choice of $\calh$) is closely 
related to the $p$-local homotopy theory of $BG$.

Fusion and linking categories were designed to a large extent to capture 
the ``$p$-local structure'' of finite groups, blocks, and $p$-completed 
classifying spaces in a way which does not depend directly on the 
structure of the ambient group. Many results in group theory, such as 
Alperin's fusion theorem \cite{alperin67} and the work by Alperin and 
Brou\'e on fusion in block theory \cite{AB79}, can be formulated in terms 
of fusion categories. One is thus led to search for an axiomatic 
definition of these concepts. The definition of a \emph{saturated fusion 
system} $\calf$ over a $p$-group $S$, generalizing $p$-fusion categories 
of finite groups, was first given by Puig \cite{Puig}. A simplified (but 
equivalent) definition of a saturated fusion system,  along with  an 
axiomatic definition of a ``centric'' linking system, was later given in 
\cite[\S1]{BLO2}. Here, the word ``centric'' refers to the set of objects 
in the linking system, which will be described in Section \ref{review}.  A 
\emph{$p$-local finite group} is then defined to be a triple 
$(S,\calf,\call)$, where $S$ is  a finite $p$-group, $\calf$ is a 
saturated fusion system over $S$, and $\call$ is a centric linking system 
associated to $\calf$.  The \emph{classifying  space} of such a triple is 
the $p$-completed nerve $|\call|\pcom$. For any $S\le G$ as above, 
$(S,\calf_S(G),\call^c_S(G))$ (where the categories are taken for 
appropriate families of subgroups) is a $p$-local finite group with 
classifying space $|\call^c_S(G)|\pcom\simeq{}BG\pcom$.  

The main goal of this paper is to examine the role of the set $\calh$ of 
subgroups of $S$ on which the fusion and linking systems are defined; 
i.e., to show when the set can be changed without changing $\calf$ and 
$\call$ in an ``essential'' way. Related questions have been studied 
extensively when $\calf$ comes from a finite group $G$, both in connection 
with the Alperin's fusion theorem (cf.\ \cite{alperin67} and 
\cite{puig76}) and more indirectly in connection with the study of 
homology decompositions (cf.\ \cite{Dwyer2} and \cite{Grodal}).  In a 
subsequent paper \cite{BCGLO2}, we use the tools developed in this paper 
to study the extension theory of fusion systems and $p$-local finite 
groups, in part motivated by our desire to develop more ways of 
constructing $p$-local finite groups that do not come from groups.  Such 
``exotic'' $p$-local finite groups do exist for all primes, and examples 
are given in \cite[\S9]{BLO2}, \cite{RV}, \cite{LO}, and \cite{BM}, but we 
still have no really good tools for constructing them, nor any sense of 
how frequently they occur.

We now describe the results of the paper in more detail. We refer the 
reader to Section \ref{review} for the definitions of abstract saturated 
fusion systems and centric linking systems; and also of $\calf$-centric 
and $\calf$-radical subgroups for a fusion system $\calf$ (analogous to 
the usual concepts of $p$-centric subgroups and radical $p$-subgroups of a 
finite group).  However, the precise definitions will not be essential to 
follow this introductory discussion. We also refer the reader to the end 
of the introduction for a list of notation which will be used throughout 
the paper.

One of the most difficult problems, when constructing exotic fusion 
systems, is showing that the fusion system one has constructed satisfies 
the axioms of saturation (see Definition \ref{sat.Frob.}).  This job is 
clearly simpler if one only needs to check the axioms on subgroups which 
are centric, rather than having to do so on all subgroups. The following 
theorem is used several times in our paper \cite{BCGLO2}, and can be used 
to shorten the proof of saturation of the exotic fusion systems in 
\cite[\S9]{BLO2}.

\begin{Th}  \label{centric->saturated-Intro}
Let $\calf$ be a fusion system over a finite $p$-group $S$, and assume
that all morphisms in $\calf$ are composites of restrictions of
morphisms between $\calf$-centric subgroups. If $\calf$ satisfies the
axioms of saturation (Definition \ref{sat.Frob.}) when applied to
$\calf$-centric subgroups of $S$, then $\calf$ is saturated.
\end{Th}

This theorem is stated more precisely, and in greater generality, in 
Theorem \ref{centric->saturated}. There, we replace ``$\calf$-centric'' 
subgroups in the above formulation by ``any collection of subgroups 
containing all subgroups which are both $\calf$-centric and 
$\calf$-radical and is closed under $\calf$-conjugacy''; but at the price 
of an additional hypothesis.  As such, it can be thought of as a converse 
to Alperin's fusion theorem for abstract fusion systems (as shown in 
\cite[2.13]{Puig} and \cite[Theorem~A.10]{BLO2}), which says that if $\calf$ 
is a saturated fusion system, then it is generated by restrictions of 
automorphisms of $\calf$-centric $\calf$-radical subgroups.  

In many applications, it is useful to construct linking systems with 
respect to different sets of subgroups than the $\calf$-centric subgroups 
of $S$. If $G$ is a finite  group, then we call a $p$-subgroup $P\le{}G$ 
\emph{$p$-quasicentric} if $O^p(C_G(P))$ has order prime to $p$; 
equivalently, if $BC_G(P)\pcom$ is the classifying space of some 
$p$-group. When $\calf$ is a  saturated fusion system over a $p$-group 
$S$, then we make an analogous definition of an 
\emph{$\calf$-quasicentric} subgroup of $S$ in Section \ref{quasicentric}. 
When $\calf$ is the fusion system of a block $b$ with defect group $S$, 
then the $\calf$-quasicentric subgroups of $S$ correspond to the 
\emph{nil-centralized} pointed groups, in the sense of Puig \cite{puig00}, 
which are associated to $b$. 

Our next theorem shows that the homotopy type of the classifying space of 
a $p$-local finite group $(S,\calf,\call)$ is also determined by a linking 
system based on \emph{any} set of $\calf$-quasicentric subgroups of $S$ 
which contains at least those which are both $\calf$-centric and 
$\calf$-radical.  This result can also be interpreted as a statement about 
homology decompositions for $p$-local finite groups, and as such is 
motivated by \cite[1.20]{Dwyer2} and \cite[Theorem 1.5]{Grodal}.  It is 
restated and proved as  Theorem \ref{|Lqc|=|L|}, and is essential when 
studying ``extensions'' of $p$-local finite groups with $p$-group quotient 
in \cite{BCGLO2}.

\begin{Th}  \label{|Lqc|=|L|-Intro}
Let $\SFL$ be a $p$-local finite group.  Then there exists a category 
$\callq$ containing $\call$ as a full subcategory, whose objects are the 
$\calf$-quasicentric subgroups of $S$, and such that the inclusion of 
nerves $|\call|\subseteq|\call^q|$ is a homotopy equivalence.  
Furthermore, if $\calh$ is any collection of $\calf$-quasicentric 
subgroups of $S$ containing all $P\le S$ which are both $\calf$-centric 
and $\calf$-radical, and $\call^\calh\subseteq\callq$ is the full 
subcategory whose objects are the subgroups in $\calh$, then the 
inclusions of $\call^\calh$ and $\call$ in $\callq$ induce homotopy 
equivalences 
	$$ |\call^\calh|\simeq|\callq|\simeq|\call|. $$
\end{Th}

We conclude this paper, in Section \ref{constrained}, with a very 
specialized family of examples: fusion systems whose entire structure is 
controlled by a single $p$-subgroup. If $G$ is a finite group which has no 
nontrivial normal subgroup of order prime to $p$, then $G$ is called 
\emph{$p$-constrained} if there is a normal $p$-subgroup $P\nsg{}G$ such 
that $C_G(P)\le P$; equivalently, such that $G/P$ can be identified (via 
conjugation) with a subgroup of $\Out(P)$.  In Section \ref{constrained}, 
we give an analogous definition of a constrained fusion system (Definition 
\ref{D:constrained}), and then prove the following proposition (restated 
as Proposition \ref{constrained1}).

\begin{Prp}\label{D:constrained-intro}
Let $\calf$ be a constrained saturated fusion system over a finite
$p$-group $S$. Then there exists a unique $p'$-reduced $p$-constrained
finite group $G$ such that $\calf = \calf_S(G)$.
\end{Prp}

For easy reference, we end the introduction with a list of notation and 
terminology which is used throughout the paper.

\begin{enumerate} \setlength{\itemindent}{-6pt}%
\renewcommand{\labelenumi}{$\bullet$}%
\item $\sylp{G}$ denotes the set of Sylow $p$-subgroups of $G$.
\item $O_p(G)$ is the maximal normal $p$-subgroup of $G$.
\item $O_{p'}(G)$ is the maximal normal subgroup of $G$ of order prime to
$p$.
\item $G$ is $p$-reduced ($p'$-reduced) if $O_p(G)=1$ (if $O_{p'}(G)=1$).
\item $O^p(G)$ the minimal normal subgroup of $G$ of $p$-power index.
\item $N_G(P,Q)=\{x\in{}G\,|\,xPx^{-1}\le{}Q\}$ \quad (for $P,Q\le{}G$).
\item $c_x$ denotes conjugation by $x$ ($g\mapsto{}xgx^{-1}$).
\item $N_S^K(P)=\{x\in{}N_S(P)\,|\,c_x\in{}K\}$ \quad (for $P\le{}S$ and 
$K\le\Aut(P)$).
\item $\Hom_G(H,K)$ (for $H,K\le{}G$) is the set of homomorphisms
from $H$ to $K$ induced by conjugation in $G$.
\item $\Aut_G(H)=\Hom_G(H,H)$, \ and \ $\Out_G(H)=\Aut_G(H)/\Inn(H)$.
\item In a fusion system $\calf$, $\homf(P,Q)=\Mor_\calf(P,Q)$,
$\isof(P,Q)=\homf(P,Q)$ if $|P|=|Q|$, $\autf(P)=\isof(P,P)$, \ and \
$\Out_\calf(P)=\Aut_\calf(P)/\Inn(P)$.
\end{enumerate}

All five authors would like to thank the University of Aberdeen, the 
University of Paris 13 and the CRM in Barcelona for giving us the 
opportunity on several occasions to get together while doing this work.


\section{A quick review of $p$-local finite groups}
\label{review}

We first recall the definitions of a fusion system, and a saturated fusion 
system, in the form given in \cite{BLO2}.

\begin{Defi}[{\cite{Puig} and \cite[Definition~1.1]{BLO2}}]
A fusion system over a finite $p$-group $S$ is a category $\calf$, where 
$\Ob(\calf)$ is the set of all subgroups of $S$, and which satisfies the 
following two properties for all $P,Q\le{}S$:
{\renewcommand{\labelenumi}{\hskip-4pt$\bullet$}
\begin{enumerate}
\item  $\Hom_S(P,Q) \subseteq \homf(P,Q) \subseteq \Inj(P,Q)$; and
\item  each $\varphi\in\homf(P,Q)$ is the composite of an isomorphism in
$\calf$ followed by an inclusion.
\end{enumerate}
}
\end{Defi}

We next specify certain collections of subgroups relative to a given 
fusion system. If $\calf$ is a fusion system over a finite $p$-subgroup 
$S$, then two subgroups $P,Q\le S$ are said to be $\calf$-conjugate if 
they are isomorphic as objects of the category $\calf$. 

\begin{Defi}\label{centric-radical-def}
Let $\calf$ be a fusion system over a finite
$p$-subgroup $S$. 
\renewcommand{\labelenumi}{\hskip-4pt$\bullet$}%
\begin{enumerate}
\item A subgroup $P \leq S$ is $\calf$-centric if $C_S(P') =Z(P')$ 
for all $P' \leq S$ which are $\calf$-conjugate to $P$. 
\item A subgroup $P \leq S$ is $\calf$-radical if $\Out_\calf(P)$ 
is $p$-reduced; i.e., if $O_p(\outf(P))=1$. 
\end{enumerate}
\end{Defi}

If $\calf = \calf_S(G)$ for some finite group $G$, then $P\le S$ is 
$\calf$-centric if and only if $P$ is $p$-centric in $G$ (i.e., 
$Z(P)\in\Syl_p(C_G(P))$), and $P$ is $\calf$-radical if and only if 
$N_G(P)/P\cdot C_G(P)$ is $p$-reduced.

The following additional definitions and conditions are needed in order for 
these systems to be very useful.

\begin{Defi}[{\cite{Puig}, see \cite[Definition~1.2]{BLO2}}] \label{sat.Frob.}
\renewcommand{\labelenumi}{\hskip-4pt$\bullet$}
Let $\calf$ be a fusion system over a $p$-group $S$.
\begin{enumerate}
\item A subgroup $P\le{}S$ is \emph{fully centralized in $\calf$} if
$|C_S(P)|\ge|C_S(P')|$ for all $P'\le{}S$ which is $\calf$-conjugate to
$P$.
\item A subgroup $P\le{}S$ is \emph{fully normalized in $\calf$} if
$|N_S(P)|\ge|N_S(P')|$ for all $P'\le{}S$ which is $\calf$-conjugate to
$P$.
\item $\calf$ is a \emph{saturated fusion system} if the following
two conditions hold:
\begin{enumtwo}  \renewcommand{\labelenumii}{\textup{(\Roman{enumii})}}
\item For all $P\le{}S$ which is fully normalized in $\calf$, $P$ is fully
centralized in $\calf$ and $\Aut_S(P)\in\sylp{\autf(P)}$.
\item If $P\le{}S$ and $\varphi\in\homf(P,S)$ are such that $\varphi{}P$ is
fully centralized, and if we set
    $$ N_\varphi = \{ g\in{}N_S(P) \,|\, \varphi c_g\varphi^{-1} \in
    \Aut_S(\varphi{}P) \}, $$
then there is $\widebar{\varphi}\in\homf(N_\varphi,S)$ such that
$\widebar{\varphi}|_P=\varphi$.
\end{enumtwo}
\end{enumerate}
\end{Defi}

If $G$ is a finite group and $S\in\sylp{G}$, then the category 
$\calf_S(G)$ defined in the introduction is a saturated fusion system (see 
\cite[Proposition~1.3]{BLO2}).

We now turn to linking systems associated to abstract fusion systems.

\begin{Defi}[{\cite[Definition~1.7]{BLO2}}]  \label{L-cat}
Let $\calf$ be a fusion system over the $p$-group $S$.  A \emph{centric
linking system associated to $\calf$} is a category $\call$ whose objects
are the $\calf$-centric subgroups of $S$, together with a functor
$\pi \:\call\Right5{}\calf^c$,
and ``distinguished'' monomorphisms $P\Right1{\delta_P}\Aut_{\call}(P)$
for each $\calf$-centric subgroup $P\le{}S$, which satisfy the following
conditions.
\begin{enumerate}
\renewcommand{\labelenumi}{\textup{(\Alph{enumi})}}%
\item  $\pi$ is the identity on objects.  For each pair of objects
$P,Q\in\call$, $Z(P)$ acts freely on $\Mor_{\call}(P,Q)$ by composition
(upon identifying $Z(P)$ with $\delta_P(Z(P))\le\Aut_{\call}(P)$), and
$\pi$ induces a bijection
	$$ \Mor_{\call}(P,Q)/Z(P) \Right5{\cong} \homf(P,Q). $$

\item  For each $\calf$-centric subgroup $P\le{}S$ and each $x\in{}P$,
$\pi(\delta_P(x))=c_x\in\Aut_{\calf}(P)$.

\item  For each $f\in\Mor_{\call}(P,Q)$ and each $x\in{}P$,
$f\circ\delta_P(x)=\delta_Q(\pi{}f(x))\circ{}f$.
\end{enumerate}
\end{Defi}

A \emph{$p$-local finite group} is defined to be a triple $\SFL$, where $S$
is a finite $p$-group, $\calf$ is a saturated fusion system over $S$, and
$\call$ is a centric linking system associated to $\calf$.  The
\emph{classifying space} of the triple $\SFL$ is the $p$-completed nerve
$|\call|\pcom$.

For any finite group $G$ with Sylow $p$-subgroup $S$, the category 
$\call^c_S(G)$ defined in the introduction is easily seen to be a centric 
linking system associated to $\calf_S(G)$.  Thus 
$(S,\calf_S(G),\call_S^c(G))$ is a $p$-local finite group, with 
classifying space $|\call^c_S(G)|\pcom\simeq{}BG\pcom$ (see 
\cite[Proposition~1.1]{BLO1}).

The following definitions are somewhat more specialized, and are 
translations to the setting of fusion systems of the concepts of a normal 
$p$-subgroup of a finite group, and of strongly and weakly closed 
subgroups.

\begin{Defi} \label{D:normal}
Let $\calf$ be a saturated fusion system over a $p$-group $S$.  Then for
any normal subgroup $Q\nsg{}S$,
\begin{enumerate}
\item $Q$ is \emph{strongly closed in $\calf$} if no element of $Q$ is
$\calf$-conjugate to an element of $S{\sminus}Q$;
\item $Q$ is \emph{weakly closed in $\calf$} if no other subgroup of $S$ is
$\calf$-conjugate to $Q$; and
\item $Q$ is \emph{normal in $\calf$} if each morphism 
$\alpha\in\homf(P,P')$ in $\calf$ extends to some
$\widebar{\alpha}\in\homf(PQ,P'Q)$ such that $\widebar{\alpha}(Q)=Q$.
\end{enumerate}
\end{Defi}

Equivalently, $Q\nsg{}S$ is normal in $\calf$ if and only if the 
normalizer fusion system $N_{\calf}(Q)$ is equal to $\calf$ as fusion 
systems over $S$ (see \cite[Definition~6.1]{BLO2}).  The next proposition, 
which is motivated by \cite[Proposition IV.2]{puig76}, gives two equivalent 
conditions for a subgroup to be normal in $\calf$.

\begin{Prop}  \label{Q<|F}
Let $\calf$ be a fusion system over $S$.  Then the following conditions on
a subgroup $Q\le{}S$ are equivalent:
\begin{enumerate}
\item $Q$ is normal in $\calf$.
\item $Q$ is strongly closed in $\calf$ and is contained in all
$\calf$-radical subgroups of $S$.
\item $Q$ is weakly closed in $\calf$ and is contained in all
$\calf$-radical subgroups of $S$.
\end{enumerate}
\end{Prop}

\begin{proof}  Assume first that $Q$ is normal in $\calf$. In particular, 
if an element  $x\in Q$ is $\calf$-conjugate to an element $y\in 
S{\sminus}Q$, then the isomorphism in $\calf$ from $\gen{x}$ to $\gen{y}$ 
extends to a morphism $Q=\gen{Q, x}\rTo \gen{Q, y}$.  But such a morphism  
clearly cannot send $Q$  to itself.  Thus $Q$ is strongly closed in 
$\calf$.  If $P\le{}S$  does not contain $Q$, then $N_{PQ}(P)/P$ is a 
nontrivial  $p$-subgroup of $\outf(P)$, which is in fact normal there. To 
see normality notice that if $\alpha \in \Aut_\calf(P)$ then $\alpha$ 
extends to $\bar \alpha \in \Aut_{\calf}(PQ)$ since $Q$ normal in $\calf$, 
so for all $x \in N_{PQ}(P)$ we have $\alpha c_x \alpha^{-1} = (\bar 
\alpha c_x \bar \alpha^{-1})_{|P} = c_{\bar \alpha(x)} \in \Aut_{PQ}(P)$. 
Hence such a subgroup $P$ cannot be $\calf$-radical. Thus, all 
$\calf$-radical subgroups of $S$ contain $Q$.  This shows 
(a) $\Rightarrow$ (b). 

Condition (b) clearly implies (c), and so it remains to show (c) 
$\Rightarrow$ (a). Assume that $Q$ is weakly closed in $\calf$, and that 
all $\calf$-radical subgroups contain $Q$.  Then by 
Alperin's fusion theorem, each morphism in $\calf$ is a composite of 
morphisms, each of which is the restriction of a morphism between 
subgroups containing $Q$, and which necessarily sends $Q$ to itself (since 
$Q$ is weakly closed).  In other words, each $\varphi\in\homf(P,P')$ 
extends to a morphism $\widebar{\varphi}\in\homf(PQ,P'Q)$ which sends $Q$ 
to itself, and hence $Q$ is normal in $\calf$.
\end{proof}


\section{Centric and radical subgroups determine saturation}
\label{saturation}

Given a fusion system which is not known to come from a group (or a 
block), it turns out to be difficult in general to show that it is 
saturated when using the definition directly.  This is one of the 
obstacles one encounters when trying to construct $p$-local finite groups 
which do not come from groups.

The main result of this section, Theorem \ref{centric->saturated},
says that it suffices to check the axioms of saturation on the
centric subgroups, in the sense that any fusion system which
satisfies these axioms for its centric subgroups generates a
saturated fusion system in a way made precise below.  In fact, our
result is stronger than that. We prove that it suffices to check
the axioms of saturation on those subgroups which are centric and
radical, and a much weaker condition on the centric subgroups
which are not radical.

Before stating the main results, we make some definitions.

\begin{Defi}
Let $\calf$ be any fusion system over a finite $p$-group $S$, and
let $\calh$ be a set of subgroups of $S$ closed under $\calf$-conjugacy.
\begin{enumerate}

\item $\calf$ is \emph{$\calh$-generated}
if every morphism in $\calf$ is a composite of restrictions of
morphisms in $\calf$ between subgroups in $\calh$.

\item $\calf$ is \emph{$\calh$-saturated} if conditions
(I) and (II) hold in $\calf$ for all subgroups $P\in\calh$.

\end{enumerate}
\end{Defi}

In terms of these definitions, Alperin's fusion theorem for
abstract fusion systems (in the form shown in
\cite[Theorem~A.10]{BLO2}) can be reformulated by saying that if
$\calf$ is a saturated fusion system over $S$, and $\calh$ is the
family of $\calf$-centric, $\calf$-radical subgroups of $S$, then
$\calf$ is $\calh$-generated.

Our main result in this section can be thought of as a converse to
this form of the fusion theorem.  In practice, it
often simplifies the task of deciding whether a fusion system is
saturated or not. As one example, the proof of \cite[Proposition
9.1]{BLO2} --- the proof that the fusion systems constructed there
are saturated --- becomes far simpler when we can use Theorem
\ref{centric->saturated}, applied with $\calh$ the set of
$\calf$-centric subgroups of $S$.

\begin{Thm}  \label{centric->saturated}
Let $\calf$ be a fusion system over a finite $p$-group $S$, and let 
$\calh$ be a set of subgroups of $S$ closed under $\calf$-conjugacy that 
contains all $\calf$-centric, $\calf$-radical subgroups of $S$. Assume 
that $\calf$ is $\calh$-generated and $\calh$-saturated, and that
\begin{enumerate}
\item[\textup{($*$)}] each $\calf$-conjugacy class of subgroups of $S$
which are $\calf$-centric but not in $\calh$ contains at least one
subgroup $P$ such that $\Out_S(P)\cap{}O_p(\outf(P))\ne1$.
\end{enumerate}
Then $\calf$ is saturated.
\end{Thm}

Note that the condition that $\mathcal{H}$ contain all $\calf$-centric, 
$\calf$-radical subgroups of $S$ is implied by ($*$); but we keep it in 
the statement for the sake of emphasis.  Condition ($*$) might at first 
sight seem artificial, but we construct an example at the end of this 
section showing that it is, in fact, necessary, and not implied by the 
other hypotheses of the theorem.

We first discuss the relation between conditions (I) and (II) in
Definition \ref{sat.Frob.}, and certain other, similar conditions
on fusion systems. We recall the definition of $N_{\varphi}$ for
any given $\varphi \in \Mor_{\calf}(P,Q)$,
$$ N_\varphi = \{x\in{}N_S(P)\,|\, \varphi c_x\varphi^{-1} \in
    \Aut_S(\varphi(P)) \}. $$

\begin{Lem}  \label{newax}
Let $\calf$ be a fusion system over a $p$-group $S$, and let
$\calh$ be a set of subgroups of $S$ closed under
$\calf$-conjugacy. Consider the following conditions on $\calf$:
\begin{enumerate}
\item[\textup{(I)$_{\calh}$}: ]  For each subgroup
$P\in\calh$ fully normalized in $\calf$, $P$ is fully centralized and
$\Aut_S(P)\in\sylp{\autf(P)}$.

\item[\textup{(I$'$)$_{\calh}$}: ]  Each $P\in\calh$ is $\calf$-conjugate
to a fully centralized subgroup $P'\in\calh$ such that
$\Aut_S(P')\in\sylp{\autf(P')}$.

\item[\textup{(II)$_{\calh}$}: ]  For each $P\in\calh$, and each
$\varphi\in\homf(P,S)$ such that $\varphi(P)$ is fully centralized
in $\calf$, $\varphi$ extends to a morphism
$\widebar{\varphi}\in\homf(N_\varphi,S)$.

\item[\textup{(IIA)$_{\calh}$}: ]  Each $\calf$-conjugacy class
$\calp\subseteq\calh$ contains a fully normalized subgroup
$\widehat{P}\in\calp$ with the following property:  for all
$P\in\calp$, there exists a morphism
$\varphi\in\homf(N_S(P),N_S(\widehat{P}))$ such that
$\varphi(P)=\widehat{P}$.

\item[\textup{(IIB)$_{\calh}$}: ]  For each fully normalized subgroup
$\widehat{P}\in\calh$ and each $\varphi\in\autf(\widehat{P})$,
there is a morphism
$\widebar{\varphi}\in\homf(N_\varphi,N_S(\widehat{P}))$ which
extends $\varphi$.
\end{enumerate}
Then
\begin{enumerate}
\item \ \textup{(I)$_{\calh}$} $\Longleftrightarrow$
\textup{(I$'$)$_{\calh}$}; \ and \
\item \ \textup{(I)$_{\calh}$ + (II)$_{\calh}$ }
$\Longrightarrow$ \textup{ (IIA)$_{\calh}$ + (IIB)$_{\calh}$ }
$\Longrightarrow$ \textup{ (II)$_{\calh}$}.
\end{enumerate}
\end{Lem}

\begin{proof} \noindent\textbf{(a) }Condition (I)$_\calh$ clearly implies 
(I$'$)$_\calh$, since every $P\le S$ is $\calf$-conju\-gate to a fully 
normalized subgroup. To see the converse, assume $P\in\calh$ is fully 
normalized. By \textup{(I$'$)$_{\calh}$} we can choose $P'\in\calh$ which 
is $\calf$-conjugate to $P$, fully centralized, and satisfies 
$\Aut_S(P')\in\sylp{\autf(P')}$. Then
	\begin{align*} 
        |\Aut_S(P')|{\cdot}|C_S(P')| &=
        |N_S(P')| \le |N_S(P)| \\
        &= |\Aut_S(P)|{\cdot}|C_S(P)| \le
        |\Aut_S(P')|{\cdot}|C_S(P')| \,: 
	\end{align*}
the first inequality holds since $P$ is fully normalized, and the
second by the assumptions on $P'$.  Thus all of these inequalities
are equalities, and so $P$ is fully centralized and
$\Aut_S(P)\in\sylp{\autf(P)}$.

\noindent\textbf{(b)} Assume (I)$_{\calh}$ and (II)$_{\calh}$
hold; we next prove that this implies (IIA)$_{\calh}$ and
(IIB)$_{\calh}$.  We first check condition (IIB)$_{\calh}$. Let
$\widehat{\varphi}=\iota_{\widehat{P}}\circ \varphi$ where
$\iota_{\widehat{P}}$ is the inclusion of $\widehat{P}$ in $S$.
Since $\widehat{P}$ is fully normalized, condition (I)$_{\calh}$
implies that $\widehat{P}$ is also fully centralized. By condition
(II)$_{\calh}$, $\widehat{\varphi}$ extends to
$\widebar{\varphi}\in\homf(N_{\widehat{\varphi}},S)$, where
	\begin{align*} 
        N_{\widehat{\varphi}} &= \{g\in{}N_S(\widehat{P})\,|\,
        \widehat{\varphi} c_g\widehat{\varphi}^{-1} \in
        \Aut_S(\widehat{\varphi}(\widehat{P})) \} \\
        &= \{g\in{}N_S(\widehat{P})\,|\, \varphi c_g \varphi^{-1}
        \in \Aut_S(\widehat{P}) \} = N_\varphi \,. 
	\end{align*}
Furthermore, $\Im(\widebar{\varphi})\le{}N_S(\widehat{P})$, since
$\widehat{P}\nsg{}N_{\widehat{\varphi}}$.

Next we check that (IIA)$_{\calh}$ holds. Fix an $\calf$-conjugacy
class $\calp\subseteq\calh$, and choose a fully normalized
subgroup $\widehat{P}\in\calp$.  Since (I)$_{\calh}$ holds,
$\widehat{P}$ is also fully centralized, and
$\Aut_S(\widehat{P})\in \sylp{\Aut_\calf(\widehat{P})}$. Thus for
any $P\in\calp$ and any $\varphi\in\isof(P,\widehat{P})$, there
exists $\chi\in\autf(P)$ such that $\varphi\chi
\Aut_S(P)\chi^{-1}\varphi^{-1}\leq \Aut_S(\widehat{P})$. Then
        $$ N_{\varphi \chi} \defeq \{g\in{}N_S(P)\,|\, \varphi \chi c_g
        \chi^{-1}\varphi^{-1} \in \Aut_S(\widehat{P}) \}
        = N_S(P) \,,$$
and hence the morphism $\varphi\chi$ extends to
$\widebar{\varphi}\in\homf(N_S(P),S)$ by (II)$_\calh$. Then
$\widebar{\varphi}(P)=\varphi\chi(P)=\widehat{P}$, and hence
$\Im(\widebar{\varphi})\le N_S(\widehat{P})$.


It remains to prove the last implication. Assume (IIA)$_\calh$ and
(IIB)$_\calh$; we must prove (II)$_{\calh}$.  Fix $P\in\calh$ and
$\varphi\in\homf(P,S)$ such that $P'\defeq\varphi(P)$ is fully
centralized in $\calf$.  Using (IIA)$_{\calh}$, choose a fully
normalized subgroup $\widehat{P}$ which is $\calf$-conjugate to
$P,P'$, and morphisms
        $$ \psi\in\homf(N_S(P),N_S(\widehat{P}))
        \qquad\textup{and}\qquad
        \psi'\in\homf(N_S(P'),N_S(\widehat{P})) $$
such that $\psi(P)=\psi'(P')=\widehat{P}$.  Set
$\widetilde{\varphi}=(\psi'|_{P'})\circ\varphi\circ(\psi|_P)^{-1}
\in\autf(\widehat{P})$.

For each $x\in{}N_\varphi$, there exists $y\in{}N_S(P')$ such that 
$\varphi{}c_x\varphi^{-1}=c_y$ as elements of $\Aut(P')$.  Then as 
automorphisms of $\widehat{P}$, 
$\widetilde{\varphi}c_{\psi(x)}\widetilde{\varphi}^{-1} = c_{\psi'(y)}$.  
This shows that $\psi(N_\varphi)\le{}N_{\widetilde{\varphi}}$.  By 
(IIB)$_{\calh}$, $\widetilde{\varphi}$ extends to a morphism 
$\widehat{\varphi}\in\homf(N_{\widetilde{\varphi}},N_S(\widehat{P}))$.

Now fix $x\in{}N_\varphi$, and let $y\in{}N_S(P')$ be such that
$\varphi{}c_x\varphi^{-1}=c_y$ as elements of $\Aut(P')$.  The
elements $\widehat{\varphi}\psi(x),\psi'(y)\in{}N_S(\widehat{P})$
induce the same conjugation action on $\widehat{P}$, and thus
differ by an element in $C_S(\widehat{P})$. Also, since $P'$ is
fully centralized, $\psi'(C_S(P'))=C_S(\widehat{P})$, and hence
        $$ \widehat{\varphi}\psi(x)\in \psi'(y){\cdot}C_S(\widehat{P}) =
        \psi'(y{\cdot}C_S(P')) \le \psi'(N_S(P')). $$
Thus $\widehat{\varphi}\psi(N_\varphi)\le\psi'(N_S(P'))$, and so
$\widehat{\varphi}\circ\psi$ factors through some morphism
$\widebar{\varphi}\in\homf(N_\varphi,N_S(P'))$ which extends
$\varphi$.  This finishes the proof of condition (II)$_{\calh}$.
\end{proof}

As an immediate consequence of Lemma \ref{newax}, we obtain the
following alternative characterization of the conditions of
saturation: a fusion system $\calf$ over $S$ is saturated if and
only if it satisfies the conditions (I$'$)$_{\calh}$,
(IIA)$_{\calh}$ and (IIB)$_{\calh}$ where $\calh$ is the set of
all subgroups $P\le S$.

\begin{Not*}
Following the notation introduced in Lemma \ref{newax} for the conditions 
stated there, we also write (--)$_Q$ or (--)$_{>Q}$ for (--)$_\calh$ when 
$\calh=\{Q\}$ or $\calh=\{P\,|\,Q\lneqq{}P\le{}S\}$, respectively. Given a 
fusion system $\calf$ over $S$, let $\cals$ be the set of all subgroups of 
$S$. For $P\le{}S$, let $\cals_{\ge{}P}\supseteq\cals_{>P}$ be the sets of 
subgroups of $S$ which contain, or strictly contain, $P$.
\end{Not*}

We will now prove two lemmas which allow us to prove Theorem 
\ref{centric->saturated} by induction on the number of $\calf$-conjugacy 
classes of subgroups of $S$ not in $\calh$. 

\begin{Lem} \label{F->NF(P)}
Let $\calf$ be a fusion system over a finite $p$-group $S$, and
let $\calh$ be a set of subgroups of $S$ closed under $\calf$-conjugacy.
Let $\calp$ be an $\calf$-conjugacy class of subgroups of $S$ which is
maximal among those not in $\calh$.  Assume $\calf$ is
$\calh$-generated and $\calh$-saturated.  Then the following hold
for any $P\in\calp$ which is fully normalized in $\calf$:
\begin{enumerate}
\item $N_\calf(P)$ is $\cals_{>P}$-saturated.
\item Each $\varphi\in\autf(P)$ is a composite of restrictions of morphisms
in $N_\calf(P)$ between subgroups strictly containing $P$.
\item $\calf$ is $(\calh\cup\calp)$-saturated if $N_\calf(P)$ is
$\cals_{\ge{}P}$-saturated.
\end{enumerate}
\end{Lem}


\begin{proof}  By a \emph{proper $\calp$-pair} will be meant a pair
$(Q,P)$, where $P\lneqq{}Q\le{}N_S(P)$ and $P\in\calp$.  Two
proper $\calp$-pairs $(Q,P)$ and $(Q',P')$ will be called
$\calf$-conjugate if there is an isomorphism
$\varphi\in\isof(Q,Q')$ such that $\varphi(P)=P'$. A proper
$\calp$-pair $(Q,P)$ will be called \emph{fully normalized} if
$|N_{N_S(P)}(Q)|\ge|N_{N_S(P')}(Q')|$ for all $(Q',P')$ in the
same $\calf$-conjugacy class.

The proof of the lemma is based on the following statements, whose proof 
will be carried out in Steps 1 to 4.

\begin{enumerate}\renewcommand{\labelenumi}{\textup{(\arabic{enumi})}}%

\item  If $(Q,P)$ is a fully normalized proper $\calp$-pair, then $Q$ is
fully centralized in $\calf$ and
    $$ \Aut_{N_S(P)}(Q) \in \sylp{\Aut_{N_\calf(P)}(Q)}. $$

\item  For each proper $\calp$-pair $(Q,P)$, and each fully normalized proper
$\calp$-pair $(Q',P')$ which is $\calf$-conjugate to $(Q,P)$,
there is some morphism 
	$$ \psi\in\homf(N_{N_S(P)}(Q),N_{N_S(P')}(Q')) $$ 
such that $\psi(P)=P'$ and $\psi(Q)=Q'$.

\item  There is a subgroup $\widehat{P}\in\calp$ which is fully 
centralized in $\calf$, and which has the property that for all 
$P\in\calp$, there is a morphism 
$\varphi\in\homf(N_S(P),N_S(\widehat{P}))$ such that 
$\varphi(P)=\widehat{P}$.  

\item Let $(Q,P)$ be a proper $\calp$-pair such that $P$ is fully
normalized in $\calf$.  If $Q$ is fully normalized in
$N_\calf(P)$, then $(Q,P)$ is fully normalized.  If $Q$ is fully
centralized in $N_\calf(P)$, then $Q$ is fully centralized in
$\calf$.
\end{enumerate}

\noindent Note that point (3) implies that $\widehat{P}$ is fully
normalized in $\calf$, and that any other $P'\in\calp$ which is
fully normalized in $\calf$ has the same properties.

Assuming points (1)--(4) have been shown, the lemma is proven as follows:

\smallskip

\noindent\textbf{(a) } We show that conditions (I) and (II) hold
in $N_\calf(P)$ for all $Q\in \cals_{>P}$. If $Q\gneqq{}P$ is
fully normalized in $N_\calf(P)$, then the proper $\calp$-pair
$(Q,P)$ is fully normalized by (4), and hence condition (I) holds
in $N_\calf(P)$ by (1).  It remains to show condition (II). Also,
by (4) again, if $P\lneqq{}Q\le{}N_S(P)$ and $Q$ is fully
centralized in $N_\calf(P)$, then it is fully centralized in
$\calf$. Hence (II) holds automatically for morphisms
$\varphi\in\Hom_{N_\calf(P)}(Q,N_S(P))$, since it holds in
$\calf$.

\smallskip

\noindent\textbf{(b) } Fix $\varphi\in\autf(P)$.  Since $\calf$ is
$\calh$-generated, there are subgroups
        $$ P=P_0,P_1,\dots,P_k=P
        \textup{\quad in $\calp$, \qquad and \qquad}
        P_i\lneqq{}Q_i\le{}S, \qquad Q_i\in\calh $$
and morphisms $\varphi_i\in\homf(Q_i,S)$ ($0\le{}i\le{}k-1$), such
that $\varphi_i(P_i)=P_{i+1}$ and
$\varphi=\varphi_{k-1}|_{P_{k-1}}\circ\cdots\circ\varphi_0|_{P_0}$.
Upon replacing each $Q_i$ by $N_{Q_i}(P_i)\gneqq{}P_i$, we can
assume that $Q_i\le{}N_S(P_i)$.  By (3), there are morphisms
$\chi_i\in\homf(N_S(P_i),N_S(P))$ for each $i$ such that
$\chi_i(P_i)=P$, where we take $\chi_0=\chi_k$ to be the identity.
Upon replacing each $\varphi_i$ by
$\chi_i\circ\varphi_i\circ\chi_i^{-1}\in\homf(\chi_i(Q_i),S)$, we can
arrange that $P_i=P$ for all $i$.  Thus $\varphi$ is a composite
of restrictions of morphisms in $N_\calf(P)$ between subgroups
strictly containing $P$.

\smallskip

\noindent\textbf{(c) } Assume that $N_\calf(P)$ is 
$\cals_{\ge{}P}$-saturated. By Lemma \ref{newax}, it is enough to check 
that Conditions (I$'$)$_{\calp}$, (IIA)$_{\calp}$, and (IIB)$_{\calp}$ are 
satisfied in $\calf$. Condition (IIA)$_{\calp}$ follows from point ($3$). 
Since $\autf(P)=\Aut_{N_\calf(P)}(P)$, it is clear that Condition 
(IIB)$_{\calp}$ holds in $\calf$. Finally, since 
$\Aut_S(P)=\Aut_{N_S(P)}(P)$, and since the properties of $\widehat{P}$ as 
described in point (3) hold for every fully normalized subgroup, 
(I)$_\calp$ also holds, and this proves that $\calf$ is 
$(\calh\cup\calp)$-saturated.

\smallskip

In order to finish the proof, it remains to prove points (1)--(4).

\smallskip

\noindent\textbf{Step 1: }  For any proper $\calp$-pair $(Q,P)$,
let $K_P\leq\Aut(Q)$ be defined by
        $$ K_P=\{\varphi\in \Aut(Q) \,|\, \varphi(P)=P\} \,.$$
If the pair $(Q,P)$ is fully normalized, then $Q$ is fully
$K_P$-normalized in $\calf$ in the sense of \cite[Definition~A.1]{BLO2}.  
Hence by \cite[Proposition~A.2(a)]{BLO2}, $Q$ is fully centralized and
    $$ \Aut_{N_S(P)}(Q) = \Aut_S(Q)\cap{}K_P \in
    \Syl_p\bigl(\autf(Q)\cap{}K_P\bigr) = \sylp{\Aut_{N_\calf(P)}(Q)}. $$
More precisely, this follows from the proof of \cite[Proposition
A.2]{BLO2}, where we need only know that $\calf$ satisfies the
axioms of saturation on subgroups containing $Q$ and its
$\calf$-conjugates.

\smallskip

\noindent\textbf{Step 2: }  Let $(Q',P')$ be any fully normalized
proper $\calp$-pair of subgroups of $S$ which is $\calf$-conjugate
to $(Q,P)$. Let $\varphi\in\isof(Q,Q')$ such that $\varphi(P)=P'$.
Since $(Q',P')$ is fully normalized, $Q'$ is fully centralized and
$\Aut_{N_S(P')}(Q')\in\sylp{\Aut_{N_\calf(P')}(Q')}$ by (1).

Since $\varphi\Aut_{N_S(P)}(Q)\varphi^{-1}$ is a $p$-subgroup of
$\Aut_{N_\calf(P')}(Q')$, there is a morphism
$\alpha\in\Aut_{N_\calf(P')}(Q')$ such that
    $$ \alpha\varphi\Aut_{N_S(P)}(Q)\varphi^{-1}\alpha^{-1} \le
    \Aut_{N_S(P')}(Q'). $$
Since $\calf$ is $\calh$-saturated, $\alpha\varphi$ extends to a
morphism $\widetilde{\alpha\varphi}\in\homf(N_{\alpha\varphi},S)$
by (II)$_Q$, where
    $$ N_{\alpha \varphi}=\{x\in{}N_S(Q)\,|\,
    \alpha \varphi c_x\varphi^{-1}\alpha^{-1}\in\Aut_S(Q')\}
    \ge N_{N_S(P)}(Q) \,.$$
Set $\psi=\widetilde{\alpha\varphi}|_{N_{N_S(P)}(Q)}\in
\homf(N_{N_S(P)}(Q),S)$. Then $\Im(\psi)\le{}N_{N_S(P')}(Q')$ by 
construction. Moreover, $\psi\vert_Q=\alpha\varphi\vert_Q$,
$\psi\vert_P=\alpha\varphi\vert_P$, and hence $\psi(P)=P'$ and
$\psi(Q)=Q'$.

\smallskip

\noindent\textbf{Step 3: }  We first show, for any $P,P'\in\calp$,
that there is a subgroup $P''\in\calp$, and morphisms
$\psi\in\homf(N_S(P),N_S(P''))$ and
$\psi'\in\homf(N_S(P'),N_S(P''))$, such that
$\psi(P)=P''=\psi'(P')$.

Let $\calt$ be the set of all sequences
        $$ \xi=\bigr(P=P_0,Q_0,\varphi_0;P_1,Q_1,\varphi_1;\dots;
        P_{k-1},Q_{k-1},\varphi_{k-1};P_k=P'\bigr) $$
such that $P_i\lneqq{}Q_i\le{}N_S(P_i)$,
$\varphi_i\in\homf(Q_i,N_S(P_{i+1}))$, and
$\varphi_i(P_i)=P_{i+1}$.  Let $\calt_r\subseteq\calt$ be the
subset of those $\xi$ for which there is no $1\le{}i\le{}k-1$ such
that $Q_i=N_S(P_i)=\varphi_{i-1}(Q_{i-1})$.  Let
$\calt\Right1{R}\calt_r$ be the ``reduction'' map, which removes
any $P_i$ such that $Q_i=N_S(P_i)=\varphi_{i-1}(Q_{i-1})$ (and
replaces $\varphi_{i-1}$ and $\varphi_i$ by their composite).

Define
        $$ I(\xi) = \{0\le{}i\le{}k-1\,|\, Q_i\lneqq{}N_S(P_i)
        \textup{ and } \varphi_i(Q_i)\lneqq{}N_S(P_{i+1}) \}. $$
If $\xi\in\calt$ and $I(\xi)\ne\emptyset$, define
        $$ \lambda(\xi) = \min_{i\in I(\xi)}[Q_i:P_i] \ge p. $$
The main observation needed to prove point (3) is that 
there exists an element $\xi\in\calt_r$ such that 
$I(\xi)=\emptyset$.  Note first that $\calt\ne\emptyset$, since $\calf$ is 
$\calh$-generated (and since $Q\gneqq{}P$ implies $N_Q(P)\gneqq{}P$). 
Hence (by the existence of the retraction functor $R$) 
$\calt_r\ne\emptyset$.

Fix an element $\xi\in\calt_r$ such that $I(\xi)\ne\emptyset$.  We
will construct $\widehat{\xi}\in\calt_r$ such that either
$I(\widehat{\xi})=\emptyset$, or
$\lambda(\widehat{\xi})>\lambda(\xi)$. For each $i\in{}I(\xi)$,
choose a fully normalized proper $\calp$-pair $(Q_i'',P_i'')$ which is
$\calf$-conjugate to $(Q_i,P_i)$, and apply (2) to choose
homomorphisms
    $$ \psi_i\in\homf\bigl(N_{N_S(P_i)}(Q_i),S\bigr)
    \qquad\textup{and}\qquad
    \psi'_i\in\homf\bigl(N_{N_S(P_{i+1})}(\varphi_i(Q_i)),S\bigr) $$
such that $\psi_i(P_i)=\psi'_i(P_{i+1})=P_i''$ and 
$\psi_i(Q_i)=\psi'_i(Q_{i+1})=Q_i''$.  Set
    $$ \widetilde{Q}_i=N_{N_S(P_i)}(Q_i)\gneqq{}Q_i
    \qquad\textup{and}\qquad
    \widetilde{Q}'_i=\psi'_i\bigl(N_{N_S(P_{i+1})}(\varphi_i(Q_i))\bigr)
    \gneqq \psi'_i(Q_i). $$
Note that if $(Q,P)$ is a proper $\calp$-pair with
$P\lneqq{}Q\lneqq{}N_S(P)$, then $N_{N_S(P)}(Q)\gneqq Q$. Thus
upon replacing the sequence $(P_i,Q_i,\varphi_i;P_{i+1})$ in $\xi$
by
    $$ (P_i,\widetilde{Q}_i,\psi_i;
    P_i'',\widetilde{Q}'_i,(\psi'_i)^{-1};P_{i+1}) $$
and similarly for the other components of $I(\xi)$, we obtain a
new element $\xi'\in\calt$, such that either
$I(\xi')=\emptyset$ or
$\lambda(\xi')>\lambda(\xi)$ (by construction
$[\widetilde{Q}_i:P_i]>[Q_i:P_i]$ and
$[\widetilde{Q}'_i:P']>[Q_i:P_i]$). Then
$\widehat{\xi}=R(\xi')\in\calt_r$ is also such that either
$I(\widehat{\xi})=\emptyset$ or
$\lambda(\widehat{\xi})>\lambda(\xi)$.

Since the function $\lambda$ is bounded above, it follows by induction that
there is $\xi \in\calt_r$ such that $I(\xi)=\emptyset$.  Write
    $$ \xi = \bigl(P_0,Q_0,\varphi_0; \dots ;
    P_{k-1},Q_{k-1},\varphi_{k-1};
    P_k\bigr) \in \calt_r
    \hskip2cm (P_0=P,\ P_k=P'). $$
The assumption $I(\xi)=\emptyset$ implies that for each $i$,
either $Q_i=N_S(P_i)$ (hence $|N_S(P_i)|\le|N_S(P_{i+1})|$), or
$\varphi_i(Q_i)=N_S(P_{i+1})$ (hence
$|N_S(P_i)|\ge|N_S(P_{i+1})|$).

Thus when $\xi \in \calt_r$, there is no $1\le{}i\le{}k-1$ such that 
$|N_S(P_i)|<|N_S(P_{i-1})|$ and also $|N_S(P_i)|<|N_S(P_{i+1})|$.  So if 
we choose $0\le{}j\le{}k$ such that $|N_S(P_j)|$ is maximal, then
    $$ |N_S(P)|\le|N_S(P_1)|\le|N_S(P_2)|\le \cdots \le |N_S(P_j)|, $$
and
    $$ |N_S(P_j)|\ge|N_S(P_{j-1})|\ge \cdots \ge|N_S(P_{k-1})|\ge
    |N_S(P')|. $$
Since $I(\xi)=\emptyset$, this implies that $Q_i=N_S(P_i)$ for all
$i<j$, and that $\varphi_i(Q_i)=N_S(P_{i+1})$ for all
$j\le{}i\le{}k-1$.  So upon setting $P''=P_j$, we obtain homomorphisms
    $$ \psi=\varphi_{j-1}\circ\cdots\circ\varphi_0 \in
    \homf(N_S(P),N_S(P'')) $$
and
    $$ \psi'=(\varphi_{k-1}\circ\cdots\circ\varphi_j)^{-1} \in
    \homf(N_S(P'),N_S(P'')) $$
such that $\psi(P)=P''=\psi'(P')$.

This was shown for an arbitrary pair of subgroups $P,P'\in\calp$. By 
successively applying the above construction to the subgroups in the 
$\calf$-conjugacy class $\calp$,  it now follows easily that there is some 
$\widehat{P}\in\calp$ such that for all $P\in\calp$, there is a morphism 
$\varphi\in\homf(N_S(P),N_S(\widehat{P}))$ such that 
$\varphi(P)=\widehat{P}$.  Note that $\widehat{P}$ is fully normalized 
since $N_S(\widehat{P})$ contains an injective image of any other $N_S(P)$ 
for $P\in \calp$. For the same reason, $\widehat{P}$ is fully centralized 
in $\calf$: its centralizer contains an injective image of the centralizer 
of any other subgroup in the $\calf$-conjugacy class $\calp$.

\smallskip

\noindent\textbf{Step 4: }  Fix a proper $\calp$-pair $(Q,P)$ such
that $P$ is fully normalized in $\calf$.  By (3), the pair
$(N_S(P),P)$ is $\calf$-conjugate to
$(N_S(\widehat{P}),\widehat{P})$.  Hence for every $P'\in\calp$,
there is
        $$ \psi\in\homf(N_S(P'),N_S(P)) $$
such that $\psi(P')=P$.

Assume $Q$ is fully normalized in $N_\calf(P)$.  Let $(Q',P')$ be
any proper $\calp$-pair $\calf$-conjugate to $(Q,P)$, and choose
$\psi$ as above.  Set $Q''=\psi(Q')$. Then $\psi$ sends
$N_{N_S(P')}(Q')$ injectively into $N_{N_S(P)}(Q'')$.  So
        $$ |N_{N_S(P')}(Q')| \le |N_{N_S(P)}(Q'')| \le |N_{N_S(P)}(Q)|; $$
where the last inequality holds since $Q$ is fully normalized in
$N_\calf(P)$.  This shows that the pair $(P,Q)$ is fully
normalized.

Finally, assume $Q$ is fully centralized in $N_\calf(P)$, and let
$Q'$ be any other subgroup in the $\calf$-conjugacy class of $Q$.
Fix $\varphi\in\isof(Q,Q')$, and set $P'=\varphi(P)$.  Again,
choose $\psi$ as above, and set $Q''=\psi(Q')$.  Then
$|C_S(Q')|\le|C_S(Q'')|$ since $\psi$ sends the first subgroup
injectively into the second, and $|C_S(Q'')|\le|C_S(Q)|$ since $Q$
is fully centralized in $N_\calf(P)$ and the pairs $(Q,P)$ and
$(Q'',P)$ are $\calf$-conjugate.  This shows that $Q$ is fully
centralized in $\calf$.
\end{proof}

Lemma \ref{F->NF(P)} reduces the problem of proving
$\calp$-saturation, for an $\calf$-conjugacy class $\calp$, to the
case where $\calp=\{P\}$ and $P$ is normal in $\calf$.  This case
is handled in the next lemma.

\begin{Lem} \label{P*}
Let $\calf$ be a fusion system over a $p$-group $S$.  Assume that
$P\nsg{}S$ is normal in $\calf$, and that $\calf$ is
$\cals_{>P}$-generated and $\cals_{>P}$-saturated.  Assume
furthermore that either $P$ is not $\calf$-centric, or
$\Out_S(P)\cap{}O_p(\outf(P))\ne1$.  Then $\calf$ is
$\cals_{\ge{}P}$-saturated.
\end{Lem}

\begin{proof}  Define
        $$ P^* = \{x\in{}S\,|\, c_x\in{}O_p(\autf(P)) \} \,. $$
It follows from the definition that $P^* \nsg S$, and we claim
that $P^*$ is strongly closed in $\calf$. Assume that $x\in P^*$
is $\calf$-conjugate to $y\in S$. Since $P$ is normal in $\calf$,
there exists $\psi \in \homf(\gen{x,P},\gen{y,P})$ which satisfies 
$\psi(P)=P$ and $\psi(x)=y$.  In particular, 
$\psi\circ c_x \circ \psi^{-1} = c_y$. It follows that $y\in P^*$, since
$c_x\in{}O_p(\autf(P))$.

Note also that $P^* \ge P{\cdot}C_S(P)$.  So by the assumption
$\Out_S(P)\cap{}O_p(\outf(P))\ne1$ if $P$ is $\calf$-centric, or
by definition if $P$ is not $\calf$-centric, $P \lneqq P^*$ in all cases.

Since $\calf$ is assumed to be $\cals_{>P}$-saturated, we need only to 
prove conditions (I)$_P$ and (II)$_P$. We first prove that these 
conditions follow from the following statement:
\begin{enumerate}  
\item[($**$)] each $\varphi\in\autf(P)$ extends to some 
$\widebar{\varphi}\in\autf(P^*)$.
\end{enumerate}
Since $P$ is normal in $\calf$, it is the only subgroup in its 
$\calf$-conjugacy class, and hence it is fully centralized and fully 
normalized. It is also clear that $P^*$ is fully normalized in $\calf$, 
since $P^*\nsg{}S$. Hence $\Aut_S(P^*)\in\sylp{\autf(P^*)}$ by (I)$_{>P}$. 
The restriction map from $\autf(P^*)$ to $\autf(P)$ is surjective by 
($**$), and so $\Aut_S(P)\in\sylp{\autf(P)}$.  Therefore condition (I)$_P$ 
holds.

Next we prove condition (II)$_P$: that each automorphism
$\varphi\in\autf(P)$ extends to a morphism defined on $N_\varphi$. By 
($**$), $\varphi$ extends to some $\psi\in\autf(P^*)$. 
Consider the groups of automorphisms
    \begin{align*}
    K &= \bigl\{\chi\in\Aut_S(P^*)\,\big|\, \chi|_P=c_x
    \textup{ some }x\in{}N_\varphi \bigr\} \\
    K_0 &= \bigl\{\chi\in\autf(P^*)\,\big|\, \chi|_P=\Id_P \bigr\} \nsg
    \autf(P^*).
    \end{align*}
By definition, for all $x\in{}N_\varphi$, we have
$(\psi{}c_x\psi^{-1})|_{P}=\chi|_{P}$ for some
$\chi\in\Aut_S(P^*)$.  In other words, as subgroups of
$\Aut(P^*)$,
    $$ \psi{}K\psi^{-1}\le \bigl\{\psi c_x\psi^{-1} \,\big|\,
    x\in{}N_\varphi \bigr\} \cdot (\psi{}K_0\psi^{-1}) \le
    \Aut_S(P^*) \cdot (\psi{}K_0\psi^{-1}). $$

In general, if $S\in\sylp{G}$, $H\nsg{}G$, and $P\le{}SH$ is a
$p$-subgroup, then there is $x\in{}H$ such that $P\le{}xSx^{-1}$.
Applied to this situation (with $G=\autf(P^*)$, $S=\Aut_S(P^*)$,
$H=\psi K_0 \psi^{-1}$, and $P=\psi K \psi^{-1}$), we see that
there is $\chi\in{}K_0$ such that
        $$ (\psi\chi)K(\psi\chi)^{-1} =
        (\psi\chi\psi^{-1}) (\psi K\psi^{-1}) (\psi\chi\psi^{-1})^{-1}
        \le \Aut_S(P^*). $$
Also, $P^*$ is fully centralized in $\calf$ by (I)$_{>P}$, since
$P^*$ is fully normalized.  So by (II)$_{>P}$,
$\psi\chi\in\autf(P^*)$ extends to a morphism $\widebar{\varphi}$
defined on $N_S^K(P^*)\ge{}N_\varphi$, and
$\widebar{\varphi}|_P=\psi|_P=\varphi$ since $\chi|_P=\Id_P$.

In order to finish the proof, it remains to prove ($**$).  Since any 
$\varphi\in\autf(P)$ is a composite of automorphisms of $P$ which extend 
to strictly larger subgroups, it suffices to show ($**$) when $\varphi$ 
itself extends to $\widetilde{\varphi}\in\isof(Q_1,Q_2)$, where 
$Q_i\gneqq{}P$.  Note that 
	\beq \widetilde{\varphi}(Q_1\cap{}P^*) = Q_2\cap{}P^* \tag{1} \eeq
since $P^*$ is strongly closed in $\calf$.  

We show ($**$) by induction on the index 
$[P^*{:}P^*\cap{}Q_1]=[P^*{:}P^*\cap{}Q_2]$.  If this index is 1, i.e., if 
$Q_1\ge{}P^*$, then $\widetilde{\varphi}(P^*)=P^*$ by (1), and hence 
$\widebar{\varphi}\defeq\widetilde{\varphi}|_{P^*}$ lies in $\autf(P^*)$ 
and extends $\varphi$.

Now assume $Q_1\ngeq{}P^*$, let $Q_3$ be any subgroup $\calf$-conjugate 
to $Q_1$ and $Q_2$ and fully normalized in $\calf$, and fix
$\varphi\in\isof(Q_2,Q_3)$.  Upon replacing $\widetilde{\varphi}$ by 
$\psi$ and by $\psi\circ\widetilde{\varphi}$, we are reduced to proving
the result when the target group is fully normalized.  So 
assume $Q_2$ is fully normalized (and hence, by (I)$_{>P}$, fully 
centralized).

This time, consider the groups of automorphisms
	\begin{align*}  
	K&=\bigl\{\chi\in\autf(Q_2)\,\big|\, \chi|_{P}\in O_p(\autf(P)) 
	\bigr\} \\
	K_0&=\bigl\{\chi\in\autf(Q_2)\,\big|\, \chi|_{P}=\Id_P \bigr\} \,.
	\end{align*}
Both $K$ and $K_0$ are normal subgroups of $\autf(Q_2)$. Also, $K/K_0$ is 
a $p$-group, since there is a monomorphism $K/K_0\rightarrow 
O_p(\autf(Q_2))$.  So any two Sylow $p$-subgroups of $K$ are conjugate by 
an element of $K_0$. 

By definition, $\Aut_{P^*}(Q_1)$ is a $p$-subgroup of $\autf(Q_1)$, all of 
whose elements restrict to elements of $O_p(\autf(P))$.  Hence 
$\widetilde{\varphi}\Aut_{P^*}(Q_1)\widetilde{\varphi}^{-1}$ is a 
$p$-subgroup of $K$.  Since $Q_2$ is fully normalized, $\Aut_S(Q_2)\in 
\sylp{\autf(Q_2)}$, and hence $\Aut_{P^*}(Q_2)=K\cap\Aut_S(Q_2)$ is a Sylow 
$p$-subgroup of $K$.  Thus there is $\chi\in{}K_0$ such that
	$$ \chi\widetilde{\varphi}\Aut_{P^*}(Q_1)
	\widetilde{\varphi}^{-1}\chi^{-1} \le \Aut_{P^*}(Q_2). $$
In particular, $N_{P^*Q_1}(Q_1)\le N_{\chi\widetilde{\varphi}}$.  Since 
$Q_2$ is fully centralized, condition (II)$_{>P}$ now implies that 
$\chi\widetilde{\varphi}$ extends to a morphism 
$\widetilde{\varphi}'\in\homf(Q'_1,N_S(Q_2))$, where 
$Q'_1=N_{P^*Q_1}(Q_1)$. Furthermore, 
$\widetilde{\varphi}'|_P=\widetilde{\varphi}|_P$ since $\chi\in{}K_0$.

By assumption, $P^*Q_1\gneqq{}Q_1$, and so 
$Q'_1=N_{P^*Q_1}(Q_1)\gneqq{}Q_1$.  Also, $Q'_1$ is generated by $Q_1$ and 
$Q'_1\cap{}P^*$ since $Q_1\le{}Q'_1\le{}P^*Q_1$.  Hence 
$Q'_1\cap{}P^*\gneqq{}Q_1\cap{}P^*$.  This shows that
	$$ [P^*{:}P^*\cap{}Q'_1] < [P^*{:}P^*\cap{}Q_1], $$
and so ($**$) now follows by the induction hypothesis.
\end{proof}

We are now ready to prove Theorem \ref{centric->saturated}.

\begin{proof}[Proof of Theorem \ref{centric->saturated}]
We are given a set $\calh$ of subsets of $S$, closed under 
$\calf$-conjugacy, such that $\calf$ is $\calh$-generated and 
$\calh$-saturated, and such that condition 
\begin{enumerate}
\item[\textup{($*$)}] each $\calf$-conjugacy class of subgroups of $S$
which are $\calf$-centric but not in $\calh$ contains at least one
subgroup $P$ such that $\Out_S(P)\cap{}O_p(\outf(P))\ne1$.
\end{enumerate}
holds.  We will prove, by induction on the number of $\calf$-conjugacy 
classes of subgroups of $S$ not in $\calh$, that $\calf$ is saturated.  If 
$\calh$ contains all subgroups, then we are done.  Otherwise, let $\calp$ 
be any $\calf$-conjugacy class of subgroups of $S$ which is maximal among 
those not in $\calh$.  We will show that $\calf$ is also 
$(\calh\cup\calp)$-saturated. Since $\calf$ is clearly 
$(\calh\cup\calp)$-generated, the result then follows by the induction 
hypothesis.

By Lemma \ref{F->NF(P)}, for any fully normalized subgroup
$P\in\calp$, the normalizer fusion system $N_\calf(P)$ is
$\cals_{>P}$-saturated, and $\autf(P)$ is generated by
restrictions of morphisms in $N_\calf(P)$ between subgroups of
$N_S(P)$ which strictly contain $P$.

Let $\calf_0$ be the fusion system over $S_0\defeq{}N_S(P)$
generated by the restriction of $N_\calf(P)$ to $\cals_{>P}$, that
is, the smallest fusion system over $S_0$ for which morphisms
between subgroups in $\cals_{>P}$ are the same as those in
$N_\calf(P)$. Then $\Aut_{\calf_0}(P)=\autf(P)$, and $\calf_0$ is
$\cals_{>P}$-saturated and $\cals_{>P}$-generated. Also, by the
assumption ($*$), either $P$ is not centric in $\calf$ (hence not
centric in $\calf_0$), or
$\Out_S(P)\cap{}O_p(\Out_{\calf_0}(P))\ne{}1$. Then $\calf_0$ is
$\cals_{\ge{}P}$-saturated by Lemma \ref{P*}, and so $\calf$ is
$(\calh\cup\calp)$-saturated by Lemma \ref{F->NF(P)} again.
\end{proof}

We end this section with a description of a example which shows
why the assumption ($*$) in Theorem \ref{centric->saturated}
($\out_S(P)\cap O_p(\outf(P))\neq 1$ if $P$ is not centric) is
needed.  We use the following standard notation:  if $k$ is a finite
field, and $n\ge1$, then $\Sigma L_n(k)$ denotes the semidirect
product of $SL_n(k)$ with the group of field automorphisms of $k$.  This 
group has an obvious action on the vector space $k^n$ and on the projective 
space $\mathbb{P}(k^n)$.  It is not
hard to see that $\Sigma{}L_2(\F_4)\cong{}S_5$: via its permutation action 
on the five points in $\mathbb{P}(\F_4{}^2)$.

Let $\Gamma=\F_4^2\rtimes S_5$, where $S_5$ acts on $\F_4^2$ via
the above isomorphism. Note that $\Gamma$ can be identified with
the subgroup of $\Sigma{}L_3(\F_4)$ generated by matrices with
bottom row $(0,0,1)$ and the field automorphism. Therefore
$\Gamma$ acts faithfully on $P=\F_4^3$.

We will define a fusion system $\calf$ over $S=P\rtimes
S'$, where $S'=\gen{(1\,2),(4\,5)} \lneqq S_5 \le \Gamma$. Consider
the following subgroups of $S$: $Q_1=P\rtimes\gen{(1\,2)}$,
$Q_2=P\rtimes\gen{(4\,5)}$, and $Q_3=P\rtimes\gen{(1\,2)(4\,5)}$.  We 
regard all of these groups, including $\Gamma$, as subgroups of 
$P\rtimes\Gamma$.

To define the morphisms in the fusion system $\calf$, let $x\in
O_2(\Gamma)\cong \F_4^2$ be the element of order two which
centralizes $S'$, and consider the subgroups $R_1=\gen{S',(3\,4\,5)}$,
$R_2=\gen{S',(1\,2\,3)}$, and $R_2'=xR_2x^{-1}$.  Set $\outf(S)=1$,
$\autf(Q_1)=\Aut_{PR_1}(Q_1)$, $\autf(Q_2)=\Aut_{PR'_2}(Q_2)$, and
$\autf(Q_3)=\Aut_S(Q_3)$. All other morphisms in the fusion system
are restrictions of the ones just described. Note in particular 
that $\outf(Q_1)\cong S_3$, $\outf(Q_2)\cong S_3$, and
$\autf(P)=\gen{R_1,R_2'}=\Gamma$. The last equality holds since
$\gen{P,R_1,R'_2}/P=\gen{S',(1\,2\,3),(3\,4\,5)}=S_5$; and
$\gen{R_1,R'_2}$ cannot be a splitting of $\Gamma/P$ in $\Gamma$
since any splitting containing $S'$ must be $P$-conjugate to the
given $S_5\le\Gamma$; so $\gen{R_1,R'_2}\cap{}P\ne1$, and
$\gen{R_1,R'_2}\ge{}P$ since $P$ is irreducible as an
$S_5$-representation.

Consider the set of subgroups $\calh = \{ S,Q_1,Q_2,Q_3 \}$. It
follows from the above description of morphisms in $\calf$ that
the subgroups in $\calh$ are the only $\calf$-centric,
$\calf$-radical subgroups. Also, $\calf$ is $\calh$-generated by
construction, and one can check that $\calf$ is $\calh$-saturated.
But $\calf$ is not saturated, since axiom (I)$_P$ fails:
$\Aut_S(P)\notin\Syl_2(\autf(P))$ since $\Aut_S(P)\cong{}C_2^2$ and
$\autf(P)\cong\Gamma$.  (One can also show that (II)$_P$ fails.)
Note that $\Out_S(P)\cap O_2(\outf(P))=S'\cap O_2(\Gamma)=1$, so
Condition ($*$) in Theorem \ref{centric->saturated} does not hold.


\section{Expanding and restricting the classifying space: 
quasicentric subgroups}
\label{quasicentric}

The goal of this section is to show how the centric linking system of a 
$p$-local finite group $\SFL$ can be extended to a larger category or 
restricted to a smaller one  without changing the homotopy type of the 
nerve of $\call$.

One motivation for doing this is a problem which frequently occurs when 
trying to construct maps between $p$-local finite groups.  A functor 
between fusion systems need not send centric subgroups to centric 
subgroups, in which case it cannot be lifted to a functor between associated 
centric linking systems. One could try to get around this by extending the 
linking systems to include all subgroups as objects.  There is in fact a 
natural extension of the linking system to a category whose objects are 
all subgroups of $S$, but in general the homotopy type of the 
$p$-completed nerve is not preserved by this extension.

We introduce here the collection of \emph{$\calf$-quasicentric} subgroups, 
which contains the centric subgroups and supports an associated linking 
system $\callq$ with properties analogous to those of the centric one. The 
important fact proved in this section is that the nerve of $\callq$ is 
homotopy equivalent to $|\call|$. Moreover, any full subcategory of 
$\callq$ whose object set contains all subgroups which are centric and 
radical also has nerve homotopy equivalent to $|\call|$.

\begin{Defi}
Let $\calf$ be a saturated fusion system over a $p$-group $S$.  A subgroup 
$P\le{}S$ is called \emph{$\calf$-quasicentric} if for each $P'$ which is 
fully centralized in $\calf$ and $\calf$-conjugate to $P$, the centralizer 
system $C_{\calf}(P')$ is the fusion system of the $p$-group $C_S(P')$.  
We let $\calfq\subseteq\calf$ denote the full subcategory whose objects 
are the $\calf$-quasicentric subgroups of $S$.
\end{Defi}

The simplest examples where $\calf$-quasicentric subgroups need not be 
$\calf$-centric are those where $\calf=\calf_S(S)$ is the fusion system of 
a $p$-group $S$: every subgroup of $S$ is $\calf$-quasicentric (while the 
trivial subgroup, at least, is not $\calf$-centric).  A more interesting 
example is given by considering a finite group $G$ with a normal subgroup 
$H$ of $p$-power index.  Fix $S\in\sylp{G}$, and set 
$S_0=S\cap{}H\in\sylp{H}$.  It is not too hard to see that each 
$\calf_{S_0}(H)$-quasicentric subgroup of $S_0$ is also 
$\calf_S(G)$-quasicentric, something which is not true for centric 
subgroups (consider, for instance, the case $A_6 \nsg \Sigma_6$). It was 
this last observation which initially led us to consider this class of 
subgroups. 

When $\calf$ is a saturated fusion system over $S$, a subgroup $P\le{}S$ 
is $\calf$-quasicentric if and only if there is no $P'$ which is 
$\calf$-conjugate to $P$, and no $Q\le{}C_S(P')\le{}S$ and 
$\alpha\in\autf(QP')$ of order prime to $p$ with $\alpha|_{P'}=\Id_{P'}$. 
Note that the set of $\calf$-quasicentric subgroups of $S$ is closed under 
$\calf$-conjugation and overgroups. 

There is also a homotopy theoretic characterization of 
$\calf$-quasicentric subgroups. If we define a map $f\:X\to Y$ to be 
quasicentric if the homotopy fibre of the map 
$f_\sharp\:\map(X,X)_{\Id_X}\Right1{}\map(X,Y)_f$ is homotopically 
discrete, then it turns out that $P\leq S$ is $\calf$-quasicentric in 
$\SFL$ if and only if the natural map $\thetaSb|_{BP}\: 
BP\Right1{}\nv{\call}\pcom$ is quasicentric.

\begin{Prop}
For any $p$-local finite group $\SFL$ and any $P\le S$, the
following are equivalent:
\begin{enumerate}
\item $P$ is $\calf$-quasicentric.
\item There is a fully centralized subgroup $P'\leq S$ which is
$\calf$-conjugate
to $P$ and such that 
	$$ \map(BP,\nv{\call}\pcom)_{\thetaS |_{BP}}\simeq
	\map(BP',\nv{\call}\pcom)_{\thetaS |_{BP'}}\simeq BC_S(P'). $$
\item The homotopy fibre of the map
$\map(BP,BP)_{\Id_{BP}}\Right1{}\map(BP,\nv{\call}\pcom)_{\thetaS |_{BP}}$
is homotopically discrete.
\item $\map(BP,\nv{\call}\pcom)_{\thetaS |_{BP}}$ is an
Eilenberg-MacLane space $K(G,1)$.
\end{enumerate}
\end{Prop}
\begin{proof} ((a)$\Rightarrow$(b)) follows by definition of
$\calf$-quasicentric and \cite[Theorem~6.3]{BLO2}.

((b)$\Rightarrow$(c)) and ((c)$\Rightarrow$(d)) follow from the long exact 
sequence of homotopy groups of the relevant fibration because 
$\map(BP,BP)_{\Id_{BP}}\simeq BZ(P)$.

Finally we prove that ((d)$\Rightarrow$(a)). Let $P'$ be a fully 
centralized subgroup of $S$ which is $\calf$-conjugate to $P$. By 
\cite[Theorem~6.3]{BLO2}, we have that 
	$$ \nv{C_{\call}(P')}\pcom\simeq 
	\map(BP',\nv{\call}\pcom)_{\thetaS |_{BP'}}\simeq 
	\map(BP,\nv{\call}\pcom)_{\thetaS |_{BP}}\simeq K(G,1). $$
In particular, $G\cong \pi_1(\nv{C_{\call}(P')}\pcom)$ is a finite 
$p$-group, and then the fusion system $C_{\calf}(P')$ coincides with the 
fusion system of $G$ (see \cite[Theorem~7.4]{BLO2}).
\end{proof}

We now turn to quasicentric linking systems; i.e., linking systems 
associated to a saturated fusion system $\calf$ whose objects are the 
$\calf$-quasicentric subgroups.  These are defined in essentially the same 
way as centric linking systems; there is just one extra axiom which is 
needed.


\begin{Defi}  \label{D:L^q}  
Let $\calf$ be a saturated fusion system over a $p$-group $S$.  A 
\emph{quasicentric linking system associated to $\calf$} consists of a 
category $\callq$ whose objects are the $\calf$-quasicentric subgroups of 
$S$, together with a functor $\pi\:\callq\rTo\calfq$, and distinguished 
monomorphisms $P{\cdot}C_S(P)\rTo^{\delta_P}\Aut_{\callq}(P)$, which 
satisfy the following conditions. 
\begin{enumerate}\renewcommand{\labelenumi}{\textup{(\Alph{enumi})$_q$}}%
\item  $\pi$ is the identity on objects and surjective on morphisms. For 
each pair of objects $P,Q$ in $\callq$ such that $P$ is fully centralized 
in $\calf$, 
$C_S(P)$ acts freely on $\Mor_{\callq}(P,Q)$ by composition (upon 
identifying $C_S(P)$ with $\delta_P(C_S(P))\le\Aut_{\callq}(P)$), and 
$\pi$ induces a bijection
    $$ \Mor_{\callq}(P,Q)/C_S(P) \Right5{\cong} \homf(P,Q). $$

\item  For each $\calf$-quasicentric subgroup $P\le{}S$ and each $g\in{}P$,
$\pi$ sends $\delta_P(g)\in\Aut_{\callq}(P)$ to
$c_g\in\Aut_{\calf}(P)$.

\item  For each $f\in\Mor_{\callq}(P,Q)$ and each $g\in{}P$, the following
square commutes in $\callq$:
    \begin{diagram}
    P & \rTo^{f} & Q \\ \dTo>{\delta_P(g)} &&
    \dTo>{\delta_Q(\pi(f)(g))} \\
    P & \rTo^{f} & Q.
    \end{diagram}

\item  For each $\calf$-quasicentric subgroup $P\le{}S$, there is some
$\iota_P\in\Mor_{\callq}(P,S)$ such that 
$\pi(\iota_P)=\incl_P^S\in\Hom(P,S)$, and such that for each 
$g\in{}P{\cdot}C_S(P)$, $\delta_S(g)\circ\iota_P=\iota_P\circ\delta_P(g)$ 
in $\Mor_{\callq}(P,S)$.
\end{enumerate}
\end{Defi}


Note that point (D)$_q$ follows from (C)$_q$ if $P$ is $\calf$-centric.  
This is why it does not appear in the definition of a centric linking 
system (Definition \ref{L-cat}).  

When $\call$ and $\callq$ are centric and quasicentric linking systems
associated to the same fusion system $\calf$, we say that $\callq$
\emph{extends} $\call$ if $\call$ is isomorphic to a full subcategory of
$\callq$ in a way which is consistent with the projection functors and the
distinguished monomorphisms.

We first show how to construct a quasicentric linking system which extends 
a given centric linking system.  One way to do this is provided by 
\cite[\S7]{BLO2}.  There, a (discrete) category $\call_{S,f}(X)$ is 
associated to any triple $(X,S,f)$, where $X$ is a space, $S$ is a 
$p$-group, and $f\:BS\Right2{}X$ is a map.  We recall this construction in 
the case where $f$ is the natural inclusion of $BS$ into 
$X=\nv{\call}\pcom$ ($f=|\theta_S|\pcom$ as defined in the next 
paragraph). As we will see, $\call_{S,f}(\nv{\call}\pcom)$ is then an 
extension of $\call$ containing all subgroups of $S$ as objects.

Let $\SFL$ be a $p$-local finite group, and let 
$\pi\:\call\Right1{}\frc{\calf}$ be the projection functor.  For each 
subgroup $P\le{}S$, let $\calb(P)$ be the category with one object $o_P$ 
and with $\textup{End}_{\calb(P)}(o_P)=P$, and identify $BP=|\calb(P)|$.  
We let $\b{g}$ denote the morphism in $\calb(P)$ corresponding to 
$g\in{}P$.  Let
    $$ \theta_P\: \calb(P) \Right4{} \call $$
be the functor which sends $o_P$ to $P$, and sends a morphism $\b{g}$ (for
$g\in{}P$) to $\delta_P(g)\in\Aut_{\call}(P)$. This induces natural maps
$\nv{\theta_P}\pcom\:BP\Right1{}\nv{\call}\pcom$. For each
$\varphi\in\Hom_{\call}(P,Q)$, we can view $\pi(\varphi)\in\homf(P,Q)$ as a 
functor $\calb(P)\to\calb(Q)$.  Let
    $$ \eta_\varphi \: \theta_P \Right3{}
    \theta_Q\circ\pi(\varphi) $$
be the natural transformation of functors given by 
	$$ \theta_P(o_P) = P \Right6{\varphi} Q = 
	\theta_Q\bigl(\pi(\varphi)(o_P)\bigr). $$
This defines an explicit
homotopy $\nv{\eta_\varphi}\:BP\times I\Right1{}\nv{\call}\pcom$ between
$\nv{\theta_P}\pcom$ and $\nv{\theta_Q}\pcom\circ B\varphi$. If for each
$\calf$-centric subgroup $P\le{}S$, we choose a morphism
$\iota_P\in\Mor_{\call}(P,S)$ which is sent to the inclusion of $P$ in $S$
by the projection functor to $\calf$, we obtain a fixed collection of natural
transformations $\eta_{\iota_P}$, and induced homotopies
$\nv{\eta_{\iota_P}}\: BP\times I\Right1{}\nv{\call}\pcom$ from
$\nv{\theta_P}\pcom$ to the restriction $\nv{\theta_S}\pcom\big |_{BP}$.

Write $f=\nv{\theta_S}\pcom$ for short. $\call_{S,\thetaS }(|\call|\pcom)$
is defined as the category  whose objects are the subgroups of $S$, and
where morphisms are
    \begin{multline*}
    \Mor_{\call_{S,\thetaS }(|\call|\pcom)}(P,Q) \\
    = \bigl\{\,(\varphi,[H])\,\bigm|\, \varphi\in\Hom(P,Q),\
    [H]\in\Mor_{\pi(\map(BP,|\call|\pcom))}
    (\thetaSb |_{BP}, \thetaSb |_{BQ}\circ{}B\varphi)\, \bigr\}. 
    \end{multline*}
Here, $\pi$ denotes the fundamental groupoid functor.  A functor
	\beq 
	\xi_{\call}\: \call \Right4{} {\call_{S,\thetaS }}(\nv{\call}\pcom)
	\eeq
is also defined as follows.  On objects, $\xi_{\call}$ is the inclusion. 
For each $\varphi\in\Mor_{\call}(P,Q)$, $ \xi_{\call}(\varphi) = 
\bigl(\pi_{P,Q}(\varphi),[H_\varphi]\bigr)$, where $H_\varphi$ is the 
homotopy $BP\times I\Right1{}\nv{\call}\pcom$ defined by
	\beq 
	H_\varphi(x,t)=\begin{cases}
	\nv{\eta_{\iota_P}}(x,1-3t) & 0\leq t\leq\frac13 \\
	\nv{\eta_\varphi}(x,3t-1) & \frac13\leq t\leq\frac23 \\
	\nv{\eta_{\iota_Q}}(B\varphi(x),3t-2) & \frac23\leq t\leq1\,.
	\end{cases}
	\eeq

By \cite[Proposition~7.3]{BLO2}, $\xi_{\call}$ defines an equivalence of 
categories to the full subcategory $\call^c_{S,\thetaS}(|\call|\pcom)
\subseteq{\call_{S,\thetaS }}(|\call|\pcom)$ whose objects
are the $\calf$-centric subgroups of $|\call|\pcom$.  In this sense, we
say that ${\call_{S,\thetaS }}(|\call|\pcom)$ is an extension of $\call$.

\begin{Prop} \label{P:L^q}  
Fix a $p$-local finite group $\SFL$, and let 
$f\:BS\Right2{}|\call|\pcom$ be as defined above.  Let
$\callq\subseteq{\call_{S,\thetaS }}(|\call|\pcom)$ be the full subcategory
whose objects are the $\calf$-quasicentric subgroups of $S$, and regard 
$\call$ as a full subcategory of $\callq$ via $\xi_\call$. Let
    $$ \pi \: \callq \Right1{} \calfq $$
be the functor which sends an $\calf$-quasicentric subgroup to itself, and
which sends a morphism $(\varphi,[H])$ to $\varphi$.  For each object $P$
in $\callq$, define the distinguished monomorphism
    $$ \delta_P\: P{\cdot}C_S(P)\Right5{} \Aut_{\callq}(P) $$
by sending $g\in{}P{\cdot}C_S(P)$ to $(c_g,[H_g])$, where $c_g$ is 
conjugation by $g$ restricted to $P$ and $H_g$ is the homotopy $BP\times 
I\Right2{\nv{\eta_g}}BS\Right2{\thetaS }\nv{\call}\pcom$ induced by the 
natural transformation $\Id\Right1{\eta_g}c_g$ which sends the unique 
object of $\calb(P)$ to the morphism $\check{g}$ of $\calb(S)$.  Then 
these structures make $\callq$ into a quasicentric linking system 
associated to $\calf$ which extends $\call$.
\end{Prop}

\begin{proof}  If $P\le{}S$ is $\calf$-quasicentric and fully centralized 
in $\calf$, then
	$$ \map(BP,|\call|\pcom)_{f|BP} \simeq |C_\call(P)|\pcom \simeq 
	BC_S(P): $$
the first equivalence by \cite[Theorem 6.3]{BLO2} and the second since $P$ 
is $\calf$-quasicentric.  Hence by definition of the category 
$\call_{S,f}(|\call|\pcom)$, $C_S(P)$ acts freely on $\Mor_{\callq}(P,Q)$ 
(for any $Q$) with orbit set $\homf(P,Q)$, and this proves axiom (A)$_q$.  
This also shows that $\xi_\call$ embeds $\call$ as a full subcategory of 
$\callq$.  

Axiom (B)$_q$ follows immediately from the definitions.  The proof of 
(C)$_q$ is identical to the argument used to prove (C) in the proof of 
\cite[Theorem~7.5]{BLO2}.  Point (D)$_q$ follows immediately from the 
construction, upon setting $\iota_P=(\incl_P^S,[c])$, where $c$ denotes 
the constant path with value $f|_{BP}\in\map(BP,|\call|\pcom)$.
\end{proof}

We are now ready to state the main result of this section:

\begin{Thm}  \label{|Lqc|=|L|}
Let $\SFL$ be a $p$-local finite group, and let $\callq$ be a quasicentric 
linking system associated to $\calf$ which extends $\call$.  Let 
$\call'\subseteq\callq$ be any full subcategory which contains all 
$\calf$-radical $\calf$-centric subgroups of $S$. Then the inclusions of 
$\call'$ and $\call$ in $\callq$ induce homotopy equivalences 
$|\call'|\simeq|\callq|\simeq|\call|$.
\end{Thm}

Theorem \ref{|Lqc|=|L|} is an immediate consequence of Proposition
\ref{|L0|simeq|L1|} below.  The rest of the section is directed towards
the proof of that proposition. We first prove some lemmas that will provide 
us with a better understanding of morphism sets in $\callq$.

\begin{Lem}  \label{L^n-prop}
Fix a $p$-local finite group $\SFL$, let $\callq$ be a quasicentric 
linking system associated to $\calf$ which extends $\call$, and let 
$\pi\:\callq\Right1{}\calfq$ be the projection.  Fix $\calf$-quasicentric 
subgroups $P,Q,R$ in $S$. Let $\varphi\in\Mor_{\callq}(P,R)$ and 
$\psi\in\Mor_{\callq}(Q,R)$ be any pair of morphisms such that 
$\Im(\pi(\varphi))\le\Im(\pi(\psi))$.  Then there is a unique morphism 
$\chi\in\Mor_{\callq}(P,Q)$ such that $\varphi=\psi\circ\chi$.
\end{Lem}

\begin{proof}  By definition of a fusion system (every morphism is the 
composite of an isomorphism followed by an inclusion), there is a unique
morphism $\widebar{\chi}\in\homf(P,Q)$ such that
$\pi(\varphi)=\pi(\psi)\circ\widebar{\chi}$.  Let
$\chi'\in\Mor_{\callq}(P,Q)$ be any morphism such that
$\pi(\chi')=\widebar{\chi}$. Choose a fully centralized group $P'$
in the $\calf$-conjugacy class of $P$ and a particular
$\alpha\in\Iso_{\callq}(P',P)$. Then by (A)$_q$, there is a unique
element $g\in{}C_S(P')$ such that
$\varphi\circ\alpha=\psi\circ\chi'\circ\alpha\circ\delta_{P'}(g)$,
and we can define
$\chi=\chi'\circ\alpha\circ\delta_{P'}(g)\circ\alpha^{-1}$.

If $\chi_1\in\Mor_{\callq}(P,Q)$ is any other morphism such that
$\varphi=\psi\circ\chi_1$, then $\pi(\chi)=\pi(\chi_1)$, hence by
(A)$_q$ again, there is a unique element $h\in{}C_S(P')$ such that
$\chi\circ\alpha=\chi_1\circ\alpha\circ\delta_{P'}(h)$; and since
$\psi\circ\chi_1\circ\alpha=
\psi\circ\chi\circ\alpha=\psi\circ\chi_1\circ\alpha\circ\delta_{P'}(h)$,
and the action of $C_S(P')$ on $\Mor_{\callq}(P',Q)$ is free, we
obtain  $h=1$ and then $\chi=\chi_1$.
\end{proof}

\begin{Lem}  \label{deltaPQ}
Fix a $p$-local finite group $\SFL$. Let $\callq$ be quasicentric linking 
system associated to $\calf$ which extends $\call$, and let 
$\pi\:\callq\Right1{}\calfq$ be the projection. Fix a choice of an 
\emph{inclusion} morphism $\iota_P\in\Mor_{\callq}(P,S)$ for each 
$\calf$-quasicentric subgroup $P\leq S$, such that $\pi(\iota_P)=\incl\in 
\Hom(P,S)$, such that the conclusion of (D)$_q$ holds if $P$ is not 
$\calf$-centric, and where $\iota_S=\Id_S$. Then, there are unique 
injections 
	$$\delta_{P,Q}\: N_S(P,Q)\Right4{} \Mor_{\callq}(P,Q)\,,$$ 
for all $\calf$-quasicentric subgroups $P,Q\leq S$, such that:
\begin{enumerate}
\item $\pi(\delta_{P,Q}(g))=c_g\in\Hom(P,Q)$, for all $g\in N_S(P,Q)$,
\item $\delta_{P,S}(1)=\iota_P$ and $\delta_{P,P}(g)=\delta_P(g)$, for all
$g\in P\cdot C_S(P)$,
\item $\delta_{Q,R}(h)\circ\delta_{P,Q}(g)=\delta_{P,R}(hg)$, for
all $g\in N_S(P,Q)$ and $h\in N_S(Q,R)$.
\end{enumerate}
\end{Lem}

\begin{proof} For each $P$ and $Q$, and each $g\in{}N_S(P,Q)$, there is by 
Lemma \ref{L^n-prop} a unique morphism $\delta_{P,Q}(g)$ such 
that 
	$$ \delta_S(g)\circ\iota_P = \iota_Q\circ\delta_{P,Q}(g). $$
We take this as the definition of the maps $\delta_{P,Q}$.  Property (a) 
follows from (B)$_q$, (b) follows by definition and the assumptions about 
$\iota_P$, and (c) follows from Lemma \ref{L^n-prop}.  (Compare with 
\cite[Proposition 1.11]{BLO2} and its proof.) 
\end{proof}

For the rest of the section, whenever we are given a $p$-local finite 
group $\SFL$, we assume that we have chosen morphisms 
$\iota_P\in\Mor_{\callq}(P,S)$, for each object $P$, such that 
$\pi(\iota_P)$ is the inclusion and the conclusion of (D)$_q$ holds.  Then 
for each $P\le{}Q$ in $\callq$, we let $\iota_P^Q\in\Mor_{\callq}(P,Q)$ be 
the unique morphism such that $\iota_P=\iota_Q\circ\iota_P^Q$ 
(Lemma~\ref{L^n-prop}).  If $\varphi\in\Mor_{\callq}(P,Q)$, and $P'\le{}P$ 
and $Q'\le{}Q$ are quasicentric subgroups such that 
$\pi(\varphi)(P')\le{}Q'$, then we write 
$\varphi|_{P'}^{Q'}\in\Mor_{\callq}(P',Q')$ for the ``restriction'' of 
$\varphi$:  the unique morphism such that 
$\iota_{Q'}^Q\circ\varphi|_{P'}^{Q'}=\varphi\circ\iota_{P'}^P$ 
(Lemma~\ref{L^n-prop} again). We also write 
$\varphi|_{P'}=\varphi|_{P'}^{Q'}$ when the target group $Q'$ is clear 
from the context.

\begin{Lem}  \label{unique-map}
Fix a saturated fusion system $\calf$ over a $p$-group $S$, and let 
$Q\le{}S$ be an $\calf$-quasicentric subgroup.  Let $P\le{}S$ be such that 
$Q\nsg{}P$, and let $\varphi,\varphi'\in\homf(P,S)$ be such that 
$\varphi|_Q=\varphi'|_Q$, and $\varphi(Q)=\varphi'(Q)$ is fully 
centralized in $\calf$.  Then there is $x\in{}C_S(\varphi(Q))$ such that 
$\varphi'=c_x\circ\varphi$.
\end{Lem}

\begin{proof} Upon replacing $P$ by $\varphi'(P)$ and $Q$ by 
$\varphi(Q)=\varphi'(Q)$, we can assume that $\varphi'=\incl_P^S$ and 
$\varphi|_Q=\Id_Q$. We are thus reduced to the case where $Q$ is fully 
centralized and $\varphi'$ is the inclusion of $P$ in $S$. 


The idea of the proof is to show that for some $x\in{}C_S(Q)$, we can extend
$\varphi\circ{}c_x$ to some $\widebar{\varphi}\in\homf(\widebar{P},S)$, for
some $\widebar{P}\ge{}P$, such that
$\widebar{\varphi}|_{\widebar{Q}}=\Id_{\widebar{Q}}$ where
$Q\lneqq\widebar{Q}\nsg\widebar{P}$.  The lemma then follows by downward
induction on $|Q|$.  Recall that the lemma holds when $Q$ is
$\calf$-centric by \cite[Proposition A.8]{BLO2}.

By definition of an $\calf$-quasicentric subgroup,
$\varphi|_{C_P(Q)}$ is conjugation by some element $x\in{}C_S(Q)$.
So after composing with $c_x$, we can assume that
$\varphi|_{C_P(Q){\cdot}Q}=\Id$.  We are thus done if
$C_P(Q){\cdot}Q\gneqq{}Q$ by taking $\widebar{P}=P$ and
$\widebar{Q}=C_P(Q){\cdot}Q$.

Assume now that $C_P(Q)\le{}Q$.  Set $K=\Aut_P(Q)$.  As in \cite[Appendix 
A]{BLO2}, we write
	$$ N_S^K(Q) = \{x\in{}N_S(Q)\,|\,c_x\in K\}, $$
and let $N_\calf^K(Q)$ be the fusion system over $N_S^K(Q)$ 
whose morphisms are defined (for $P,P'\le{}N_S^K(Q)$) by
	\begin{multline*} 
	\Hom_{N_\calf^K(Q)}(P,P') \\
	=  \bigl\{\varphi\in\homf(P,P') 
	\,\big|\, \psi|_P=\varphi,\ \psi|_Q\in{}K, \textup{ some }
	\psi\in\homf(PQ,P'Q) \bigr\} \,. 
	\end{multline*}
Then $P$, $\varphi(P)$, and $C_S(Q)$ are all 
contained in $N_S^K(Q)$. If $Q$ is not fully $K$-normalized in $\calf$, 
then there is some $\psi\in\homf(N_S^K(Q),S)$ such that $\psi(Q)$ is fully 
$\psi{}K\psi^{-1}$-normalized in $\calf$ (see \cite[Proposition 
A.2(b)]{BLO2}); and upon replacing all of these subgroups by their images 
under $\psi$, we are reduced to the case where $Q$ is fully $K$-normalized 
in $\calf$. The fusion system $N_{\calf}^K(Q)$ is saturated by 
\cite[Proposition A.6]{BLO2}; and upon replacing $\calf$ by 
$N_{\calf}^K(Q)$ we can assume that $S=N_S^K(Q)=P{\cdot}C_S(Q)$ and 
$\calf=N_{\calf}^K(Q)$.  In particular, each $\alpha\in\homf(R,R')$ 
extends to a morphism in $\homf(RQ,R'Q)$ whose restriction to $Q$ is 
conjugation by some element of $P$.

Fix $\psi\in\homf(P,S)$ such that $\psi(P)$ is fully normalized in $\calf$.
Since $\psi|_Q$ is conjugation by an element $g\in{}P$, we can replace
$\psi$ by $\psi\circ{}c_g^{-1}$, and thus arrange that $\psi|_Q=\Id$.  If
$\psi$ and $\psi\circ\varphi^{-1}$ are both conjugation by some element of
$C_S(Q)$, then so is $\varphi$; so it suffices to prove the result under
the assumption that $\varphi(P)$ is fully normalized in $\calf$.

Now, $(C_S(Q){\cdot}Q)/Q$ is a nontrivial normal subgroup of
$N_S(Q)/Q=S/Q$.  So there is an element $x\in{}C_S(Q){\sminus}Q$
such that $1\ne{}xQ\in{}Z(S/Q)$.  Then $x\in{}N_S(P)$, and acts
via the identity on $Q$ and on $P/Q$.  Thus
    $$ c_x \in \Ker\bigl[\autf(P) \Right3{} \autf(Q) \times \Aut(P/Q)
    \bigr], $$
a normal $p$-subgroup of $\autf(P)$ (see \cite[Corollary
5.3.3]{Gorenstein}). Also,
$\Aut_S(\varphi(P))\in\sylp{\autf(\varphi(P))}$ since $\varphi(P)$
is fully normalized.  Hence
$\varphi{}c_x\varphi^{-1}\in\Aut_S(\varphi(P))$ (after replacing
$\varphi$ by $\varphi \circ \xi$ where $\xi \in \autf(\varphi(P))$
if necessary). Thus, $x\in{}N_\varphi$ and $Q\lneqq N_{\varphi}$.
By (II), $\varphi$ extends to
$\widebar{\varphi}\in\homf(N_\varphi,S)$.
Now set $\widebar{P}=N_\varphi\cap{}N_S(Q)=N_\varphi$ and
$\widebar{Q}=C_{\widebar{P}}(Q){\cdot}Q$.

By construction, $x\in\widebar{Q}{\sminus}Q$. Since $Q$ is 
$\calf$-quasicentric, $\widebar{\varphi}|_{C_{\widebar{P}}(Q)}$ is 
conjugation by some element $g\in{}C_S(Q)$.  So we can replace 
$\widebar{\varphi}$ by $\widebar{\varphi}\circ{}(c_g)^{-1}$, and thus 
arrange that $\widebar{\varphi}|_{\widebar{Q}}=\Id_{\widebar{Q}}$.  Since 
$\widebar{Q}\gneqq{}Q$ and $\widebar{Q}\nsg\widebar{P}$, this finishes the 
induction step.
\end{proof}

The next lemma can be thought of as a ``lifting'' of the last one to
quasicentric linking systems.  It says that all inclusions in $\callq$ are
epimorphisms in the categorical sense.

\begin{Lem}  \label{L-epi}
Fix a $p$-local finite group $\SFL$, and let $\callq$ be a quasicentric 
linking system associated to $\calf$ which extends $\call$.  Assume 
$Q\le{}P\le{}S$ and $R\le{}S$ are $\calf$-quasicentric, and let 
$\varphi,\varphi'\in\Mor_{\callq}(P,R)$ be two morphisms such that 
$\varphi\circ\iota_Q^P=\varphi'\circ\iota_Q^P$. Then $\varphi=\varphi'$.
\end{Lem}

\begin{proof}  Since there is always a subnormal series $Q=Q_0\nsg
Q_1\nsg\cdots\nsg Q_k=P$, it suffices to prove the lemma when $Q$ is
normal in $P$.  So we assume this from now on.

It will be convenient, throughout the proof, to write
$\widehat{\alpha}=\pi(\alpha)\in\Mor(\calf)$ for any
$\alpha\in\Mor(\callq)$.  By Lemma \ref{L^n-prop},
$\varphi=\varphi'$ if and only if
$\iota_R^S\circ\varphi=\iota_R^S\circ\varphi'
\in\Mor_{\callq}(P,S)$, and similarly replacing $\varphi$ (resp.
$\varphi'$) by $\varphi \circ \iota_Q^P$ (resp. $\varphi' \circ
\iota_Q^P$).  We can thus replace $R$ by any other subgroup of $S$
which contains the images of $\widehat{\varphi}$ and
$\widehat{\varphi}'$, and in particular assume that $R\le
N_S(\widehat{\varphi}(Q))$.

The proof itself will be divided in two steps: the first dealing with a
restricted case, and the second reducing the general case to that in Step 1.

\smallskip

\noindent\textbf{Step 1:} Assume first that
$Q=\widehat{\varphi}(Q)$ and is fully normalized, and that $P$ is
fully centralized. Set
$\varphi_0=\varphi\circ\iota_Q^P=\varphi'\circ\iota_Q^P$.  By
condition (II) in Definition \ref{sat.Frob.} (and since
$Q=\widehat{\varphi}_0(Q)$ is fully centralized), there is
$\psi\in\Hom_\call(P{\cdot}C_S(Q),S)$ such that
$\widehat{\psi}|_Q=\widehat{\varphi}_0$.  Set 
$\varphi''=\psi|_P\in\Mor_{\callq}(P,S)$.  Thus
$\widehat{\varphi}''|_Q=\widehat{\varphi}_0$, so there is a unique
element $a\in{}C_S(Q)$ such that
    $$ \psi|_Q = \varphi''|_Q = 
    \iota_R^S\circ\varphi_0 \circ \delta(a). $$
We will show that $\varphi=\varphi'$ by comparing both to $\psi$ and 
$\varphi''$; the advantage of this is that condition (C) can be applied more 
easily to $\psi$.  

By Lemma \ref{unique-map}, there is some $x\in{}C_S(Q)$ 
($\widehat{\varphi}(Q)=Q$) such that $c_x\circ\widehat{\varphi} 
=\widehat{\varphi}''$. Since $P$ is fully centralized, by condition (A)$_q$ 
in Definition \ref{D:L^q} there is some $y\in{}C_S(P)$ such that
	$$ \delta_{R,S}(x)\circ\varphi = 
	\iota_R^S\circ\varphi''\circ\delta_P(y) 
	= \psi\circ\delta_{P{\cdot}C_S(Q)}(y)|_P
	= \delta_S(\widehat{\psi}(y))\circ\psi|_P
	= \delta_{S}(\widehat{\psi}(y))\circ\varphi''. $$
It follows that $\varphi''=\delta_{R,S}(z)\circ\varphi$, where
$z=\widehat{\psi}(y)^{-1}{\cdot}x\in{}C_S(Q)$.  Hence
	\begin{multline*} 
	\delta_{R,S}(\widehat{\psi}(a)^{-1}{\cdot}z)\circ\varphi_0
	= \delta_S(\widehat{\psi}(a))^{-1}\circ\varphi''|_Q
	= \delta_S(\widehat{\psi}(a))^{-1}\circ\psi|_Q \\
	= \psi\circ\delta_{P{\cdot}C_S(Q)}(a)^{-1}|_Q
	= \psi|_Q\circ\delta_Q(a)^{-1}
	= \varphi''\circ\iota_Q^P\circ\delta_Q(a)^{-1} 
	= \iota_R^S\circ \varphi_0.
	\end{multline*}
Since $\varphi_0=\iota_Q^R\circ\omega$ for some 
$\omega\in\Aut_{\callq}(Q)$, upon composing with $\omega^{-1}$,
this shows that $\delta_{Q,S}(\widehat{\psi}(a)^{-1}{\cdot}z) 
=\iota_Q^S$, and hence that $z=\widehat{\psi}(a)$.

After making a similar argument involving $\varphi'$, we now have
    $$ \delta(\widehat{\psi}(a))\circ\varphi = \varphi'' =
    \delta(\widehat{\psi}(a))\circ\varphi', $$
and this shows that $\varphi=\varphi'$.

\smallskip

\noindent\textbf{Step 2:} (General case.)
We first reduce the problem to the case in which $P$ is fully
centralized. We choose an isomorphism $\xi \in
\Mor_{\callq}(P,P')$ such that $\widehat{\xi}(P)=P'$ is fully
centralized. Upon replacing $P$ by $P'$, $\varphi$ by $\varphi
\circ \xi^{-1}$, and $\varphi'$ by $\varphi' \circ \xi^{-1}$ we
are now reduced to the case where $P$ is fully centralized in
$\calf$.

Set $Q'=\widehat{\varphi}(Q)=\widehat{\varphi}'(Q)$ for short; we
now reduce the problem to the case in which $Q=Q'$ and is fully
normalized. Let $Q''$ be any fully normalized subgroup in the
$\calf$-conjugacy class of $Q$ (and of $Q'$).  By Lemma \ref{newax}
(condition (IIB) holds), there are morphisms
    $$ \beta\in\Mor_{\callq}(N_S(Q),N_S(Q'')) \qquad\textup{and}\qquad
    \beta'\in\Mor_{\callq}(N_S(Q'),N_S(Q'')) $$
such that $\widehat{\beta}(Q)=\widehat{\beta}'(Q')=Q''$.  Set
$P''=\widehat{\beta}(P)$, and let $\beta_0\in\Iso_{\callq}(P,P'')$
be the restriction of $\beta$ (i.e., by Lemma \ref{L^n-prop} the
unique morphism such that
$\iota_{P''}^{N_S(Q'')}\circ\beta_0=\beta\circ\iota_P^{N_S(Q)}$).
Set
    $$ \psi=\beta'\circ\iota_R^{N_S(Q'')}\circ\varphi\circ\beta_0^{-1}
    \,,\
    \psi'=\beta'\circ\iota_R^{N_S(Q'')}\circ\varphi'\circ\beta_0^{-1}
    \in\Mor_{\callq}(P'',N_S(Q'')). $$
Then $\psi=\psi'$ if and only if $\varphi=\varphi'$, and
$\psi\circ\iota_{Q''}^{P''}=\psi'\circ\iota_{Q''}^{P''}$ if and
only if $\varphi\circ\iota_Q^P=\varphi'\circ\iota_Q^P$. Note that
$P''$ is $\calf$-conjugate to $P$ and the following inequality
holds:
$$|C_S(P)|=|C_{N_S(Q)}(P)|\le|C_{N_S(Q'')}(P'')|=|C_S(P'')|.$$
Since $P$ is fully centralized, it follows that
$|C_S(P)|=|C_S(P'')|$ and $P''$ is also fully centralized.

Thus, upon replacing $(Q,P,R)$ by $(Q'',P'',N_S(Q''))$, $\varphi$
by $\psi$, and $\varphi'$ by $\psi'$, we are reduced to the case
where $Q=\varphi(Q)$ is fully normalized and $P$ is fully
centralized.
\end{proof}

An immediate consequence of Lemmas \ref{L^n-prop} and \ref{L-epi} is:

\begin{Cor} \label{OL-epi}
Let $\callq$ be a quasicentric linking system associated to a saturated 
fusion system $\calf$ over a $p$-group $S$.  Then all morphisms in 
$\callq$ are monomorphisms and epimorphisms in the categorical sense.
\end{Cor}

\begin{proof} By the uniqueness in Lemma \ref{L^n-prop}, 
$\psi\circ\chi=\psi\circ\chi'$ in $\callq$ implies $\chi=\chi'$.  Hence 
all morphisms in $\callq$ are monomorphisms.  

Since each morphism in $\callq$ is the composite of an isomorphism 
followed by an inclusion, it suffices to prove that inclusions $\iota_Q^P$ 
are epimorphisms, and it clearly suffices to do this when $Q\nsg{}P$.  So 
assume $P'\le{}S$ and $\varphi,\varphi'\in\Mor_{\callq}(P,R)$ are such 
that $\varphi\circ\iota_Q^P=\varphi'\circ\iota_Q^P$.  Then 
$\iota_{P'}^S\circ\varphi=\iota_{P'}^S\circ\varphi'$ by Lemma \ref{L-epi}, 
and so $\varphi=\varphi'$ by Lemma \ref{L^n-prop}.
\end{proof}

We are now ready to prove the following proposition, of which
Theorem \ref{|Lqc|=|L|} is an immediate consequence.

\begin{Prop} \label{|L0|simeq|L1|}
Let $\SFL$ be a $p$-local finite group, and let $\callq$ be a quasicentric 
linking system associated to $\calf$ which extends $\call$.  Let 
$\call_0\subseteq\callq$ be any full subcategory such that $\Ob(\call_0)$ 
is closed under $\calf$-conjugacy.  Let $P\in\Ob(\callq)$ be maximal among 
those $\calf$-quasicentric subgroups not in $\call_0$, and let 
$\call_1\subseteq\callq$ be the full subcategory whose objects are the 
objects in $\call_0$ together with all subgroups $\calf$-conjugate to $P$. 
Assume furthermore that $P$ is not $\calf$-centric or not $\calf$-radical. 
Then the inclusion of nerves $|\call_0|\subseteq|\call_1|$ is a homotopy 
equivalence.
\end{Prop}

\newcommand{\rr}{r}
\newcommand{\callzn}{\call_0\cap N_{\call^q}(P)}
\newcommand{\callun}{\call_1\cap N_{\call^q}(P)}

\begin{proof}
Throughout the following proof, when working in any linking
system, we assume that inclusion morphisms $\iota_P^Q$ have been
chosen as in Lemma \ref{deltaPQ}.   By ``extensions'' and
``restrictions'' of morphisms we mean with respect to these
inclusions.  Also, for $\varphi\in\Mor_{\callq}(Q,Q')$, we write
$\Im(\varphi)=\Im(\pi(\varphi))\le{}Q'$ and
$\varphi(R)=\pi(\varphi)(R)\le{}Q'$ if $R\le{}Q$.

We must show that the inclusion functor
$\iota\:\call_0\rightarrow \call_1$ induces a homotopy equivalence 
$|\call_0|\simeq |\call_1|$. By Quillen's Theorem A (see \cite{Qu2}), it 
will be enough to prove that the undercategory $Q{\downarrow}\iota$ is
contractible (i.e., $|Q{\downarrow} \iota|\simeq *$) for each $Q$ in
$\call_1$. This is clear when $Q$ is not isomorphic to $P$ (since
$Q{\downarrow}\iota$ has initial object $(Q,\Id)$ in that case), so it 
suffices to consider the case $Q=P$.  Since $P$ was arbitrarily chosen in 
its isomorphism class, we can also assume that $P$ is fully normalized.

Let 
	$$ \iota_N\: \call_0\cap N_{\call^q}(P) \Right5{} 
	\call_1\cap N_{\call^q}(P) $$
be the restriction of $\iota$.  Consider the functor 
$i:P{\downarrow}\iota_N \rightarrow P{\downarrow} \iota$ induced by the 
inclusions $\call_i\cap N_{\call^q}(P)\rightarrow \call_i$ for $i=0,1$. We 
will first show that $|P{\downarrow} \iota|\simeq |P{\downarrow} \iota_N|$ 
and then that $ |P{\downarrow} \iota_N|\simeq*$. 

To prove the first statement, we construct a retraction functor
$\rr\:P{\downarrow}\iota\rightarrow P{\downarrow}\iota_N$ such that 
$\rr\circ i=\Id_{P{\downarrow}\iota_N}$, together with a natural
transformation $\bigl(i\circ\rr\rTo^\eta\Id_{P{\downarrow}i}\bigr)$.  
By Lemma \ref{newax} (condition (IIB)), for each $P'\le{}S$ which
is $\calf$-conjugate to $P$, there is a morphism in $\calf$ from
$N_S(P')$ to $N_S(P)$ which sends $P'$ isomorphically to $P$.
Hence upon lifting this to the linking system, we can choose a morphism
    $$ \Phi_{P'}\in\Mor_{\callq}(N_S(P'),N_S(P)) $$
for each such $P'$ which restricts to an isomorphism from $P'$ to
$P$. In particular, we set $\Phi_P=\Id_{N_S(P)}$.

For each nonisomorphism $\varphi\in\Mor_{\callq}(P,Q)$, set
$\widehat{\rr}(\varphi)=\Phi_{\varphi(P)}(N_Q(\varphi(P)))\gneqq P$.  We 
can factor $\varphi$ as $\varphi=\eta(\varphi)\circ\rr(\varphi)$, where
    $$ \rr(\varphi)=
    \iota_P^{\widehat{\rr}(\varphi)}\circ
    (\Phi_{\varphi(P)}|_{\varphi(P)}\circ\varphi)
    \in\Mor_{N_{\call^q}(P)}(P,\widehat{\rr}(\varphi)) $$
and
    $$ \eta(\varphi)=
    \iota_{\widebar{Q}}^Q\circ (\Phi_{\varphi(P)}|_{\widebar{Q}})^{-1} \in
    \Mor_{\call_0}(\widehat{\rr}(\varphi),Q), $$ 
where $\widebar{Q}=N_Q(\varphi(P))$. We define the functor 
$\rr:P{\downarrow} \iota \rightarrow P{\downarrow} \iota_N$ on objects by 
setting
	$$ \rr\bigl(P\Right3{\varphi}Q\bigl) = 
	\bigl(P\Right3{\rr(\varphi)}\widehat{\rr}(\varphi)\bigr). $$
For any morphism $\beta\in\Mor_{P{\downarrow}
\call_0}((Q,\varphi),(Q',\varphi'))$; i.e., for any commutative
square of the form
    \beq \begin{diagram}[w=30pt]
    P & \rTo^{\varphi} & Q \\
    \dTo<{\Id} && \dTo<{\beta} \\
    P & \rTo^{\varphi'} & Q',
    \end{diagram} \tag{1} \eeq
we claim there is a unique morphism $\widehat{\rr}(\beta)$ such that the two
squares in the following diagram commute:
    \beq \begin{diagram}[w=40pt]
    P & \rTo^{\rr(\varphi)} & \widehat{\rr}(\varphi) &
    \rTo^{\eta(\varphi)} & Q \\
    \dTo<{\Id} && \dTo<{\widehat{\rr}(\beta)} && \dTo<{\beta} \\
    P & \rTo^{\rr(\varphi')} & \widehat{\rr}(\varphi') &
    \rTo^{\eta(\varphi')} & Q' \rlap{\,.}
    \end{diagram} \tag{2} \eeq
To see this, note that by commutativity of the square (1), $\beta$
sends $N_Q(\varphi(P))$ into $N_{Q'}(\varphi'(P))$.  Hence upon
defining
    $$ \widehat{\rr}(\beta) \defeq \Phi_{\varphi'(P)} \circ \beta \circ
    \Phi_{\varphi(P)}{}^{-1} ,$$
where the three morphisms are replaced by appropriate restrictions, we get 
$\widehat{\rr}(\beta)$ such that the right square in (2) commutes. Since 
the combination of the two squares commutes by assumption, we obtain that 
$\eta(\varphi')\circ \widehat{\rr}(\beta)\circ 
\rr(\varphi)=\eta(\varphi')\circ \rr(\varphi')$, and therefore 
$\widehat{\rr}(\beta)\circ \rr(\varphi)=\rr(\varphi')$ by Lemma 
\ref{L^n-prop}. By the uniqueness of $\widehat{\rr}(\beta)$, it follows 
that this construction defines a functor, as well as a natural 
transformation $i\circ\rr\rTo^\eta\Id_{P{\downarrow}i}$.  Since 
$\rr\circ{}i=\Id_{P{\downarrow}i_N}$, this finishes the proof that 
$|P{\downarrow} \iota|\simeq |P{\downarrow} \iota_N|$. 

It remains to prove that $|P{\downarrow}\iota_N|\simeq*$. Set 
	$$ \widehat{P}= \{x\in N_S(P) \,|\, c_x\in O_p(\autf(P)) \}. $$ 
Note that $\widehat{P}\ge P{\cdot}C_S(P)$, and hence $\widehat{P}\gneqq P$ 
if $P$ is not centric. Moreover, $\widehat{P}\gneqq P$ if $P$ is not 
radical, and thus $\widehat{P}\in \call_0$ in both cases covered by the 
hypotheses of the proposition. Since $P$ is normal in $\widehat{P}$, this last
is an object in $\callzn$.

Recall that $\iota_N\:\callzn\rightarrow \callun$ denotes the inclusion. 
Let $i$ be the functor $i:\widehat{P}{\downarrow} \iota_N \rightarrow 
P{\downarrow} \iota_N$ which is induced by precomposing with the inclusion 
$\iota_{P}^{\widehat{P}}\in\Mor_{\callq}(P,\widehat{P})$. We show that $i$ 
induces a homotopy equivalence $|P{\downarrow} \iota_N|\simeq 
|\widehat{P}{\downarrow} \iota_N|$, by defining a functor 
$r\:P{\downarrow} \iota_N \rightarrow \widehat{P}{\downarrow} \iota_N$ 
such that $r\circ i=\Id_{\widehat{P}{\downarrow} \iota}$, and such that 
$i\circ r\simeq \Id_{P{\downarrow}\iota_N}$ (such that there is a natural 
transformation of functors from the identity to $i\circ r$). Then 
$|P{\downarrow} \iota_N|\simeq |\widehat{P}{\downarrow} \iota_N|$, and the 
last space is contractible since $\widehat{P}\in\callzn$. This will finish 
the proof.

Fix subgroups $Q,Q'\le N_S(P)$ containing $P$, and let
$\varphi\in\Mor_{N_{\callq}(P)}(Q,Q')$ be a morphism. 
Set $\alpha=\pi(\varphi)|_P\in\autf(P)$ for short.  Since $P$ is fully 
normalized, $\Aut_S(P)\in\sylp{\autf(P)}$, and hence 
$O_p(\autf(P))\le\Aut_S(P)$.  It follows that
	$$ N_\alpha \defeq \bigl\{x\in{}N_S(P) \,\big|\, 
	\alpha c_x\alpha^{-1}\in\Aut_S(P) \bigr\} \ge \widehat{P}; $$
and $N_\alpha\ge{}Q$ since $\alpha$ extends to 
$\pi(\varphi)\in\homf(Q,Q')$.  Thus,
since $P$ is fully centralized, $\alpha$ extends to some
$\varphi'\in\homf(Q\widehat{P},Q'\widehat{P})$ by condition (II) in 
Definition \ref{sat.Frob.}. After possibly composing 
this extension with $\delta_{Q\widehat{P}}(x)$ for some element $x\in 
C_S(P)\leq Q\widehat{P}$, we get a lifting 
$\widehat{\varphi}\in\Mor_{\callq}(Q\widehat{P}, Q'\widehat{P})$ such that 
the following diagram commutes in $\callq$:
	\begin{diagram}[w=35pt]
	P & \rTo^{\iota_P^Q} & Q & \rTo^{\varphi} & Q' \\
	\dTo<{\iota_P^{Q\widehat{P}}} && && 
	\dTo>{\iota_{Q'}^{Q'\widehat{P}}} \\
	Q\widehat{P} & \rTo^{\widehat{\varphi}} &&& Q'\widehat{P} \rlap{\,.}
	\end{diagram}
Hence by Lemma~\ref{L-epi}, $\widehat{\varphi}\circ 
\iota_{Q}^{Q\widehat{P}}= \iota_{Q'}^{Q'\widehat{P}}\circ \varphi$. This 
lifting is unique by Corollary \ref{OL-epi}; and it lies in $\callzn$, or 
in $\callun$ if $Q=P$.

The functor $r$ is defined on objects by setting
	$$ r\bigl(P\Right3{\varphi}Q\bigr)=
	\bigl(\widehat{P}\Right3{\widehat{\varphi}}Q\widehat{P}\bigr).$$ 
If $\beta:Q\rightarrow Q'$ is a morphism such that 
$\beta\circ\varphi=\varphi'$, then we define $r(\beta)=\widehat{\beta}$. 
Because of the uniqueness of the extension $\widehat{\beta}$, this 
construction defines a functor. Moreover, $r\circ 
i=\Id_{\widehat{P}{\downarrow} \iota_N}$, and $i\circ r\simeq 
\Id_{P{\downarrow}\iota_N}$, where the homotopy is induced by the natural 
transformation given by the inclusions $\iota_{Q}^{Q\widehat{P}}$.
\end{proof}

As noted above, Theorem \ref{|Lqc|=|L|} follows immediately from 
Proposition \ref{|L0|simeq|L1|}.  Another consequence of this proposition 
is the uniqueness of quasicentric linking systems associated to a given 
fusion system which extend a given centric linking system.  In fact, in 
the following proposition, we prove a slightly stronger result, by 
comparing a more general ``partial quasicentric linking system'' defined 
on smaller sets of objects with a quasicentric linking system as 
constructed in Proposition \ref{P:L^q}, and show that the first is 
contained in the second if they agree after restricting to centric radical 
subgroups.  For any $p$-local finite group $\SFL$, $\call^r$ denotes the 
full subcategory of $\call$ whose objects are the subgroups which are 
$\calf$-radical as well as $\calf$-centric.  

\begin{Prop}
Fix a $p$-local finite group $\SFL$, and let $\callq$ be the quasicentric 
linking system constructed in Proposition \ref{P:L^q}.  Let $\calh$ be any 
set of $\calf$-quasicentric subgroups of $S$ which is closed under 
$\calf$-conjugacy and overgroups, and which contains all $\calf$-centric 
$\calf$-radical subgroups.  Let $\call'$ be a category with 
$\Ob(\call')=\calh$, together with a functor $\pi'\:\call'\Right2{}\calf$, 
and distinguished monomorphisms $\delta'_P$ for all $P\in\calh$, which 
satisfy axioms (A)$_q$, (B)$_q$, (C)$_q$, and (D)$_q$ in Definition 
\ref{D:L^q}.  Assume $\call'$ contains a full subcategory isomorphic to 
$\call^r$ in a way compatible with the projection functors and 
distinguished monomorphisms.  Then $\call'$ is isomorphic to the full 
subcategory of $\callq$ with object set $\calh$, via an inclusion functor 
$\call'\Right2{}\callq$ which commutes with the projection functors and 
distinguished monomorphisms for both categories.
\end{Prop}

\begin{proof}  By Theorem \ref{|Lqc|=|L|}, 
$|\call'|\simeq|\call^r|\simeq|\call|$.  More precisely, the second 
equivalence follows directly from the theorem, and the first equivalence 
follows from the same argument applied to $\call'$, since we never needed 
to know that the linking system was defined on all $\calf$-quasicentric 
subgroups (or even on all $\calf$-centric subgroups).  

In particular, $\callq$ is a full subcategory of
	$$ \call_{S,f}(|\call|\pcom) \cong \call_{S,f}(|\call'|\pcom). $$
So if we let $\xi_{\call'}\:\call'\Right2{}\call_{S,f}(|\call'|\pcom)$
be the functor defined earlier in the section (just before Proposition 
\ref{P:L^q}), then $\xi_{\call'}$ defines an inclusion of $\call'$ into 
$\callq$, which is clearly compatible with the projection functors and 
distinguished monomorphisms.
\end{proof}


\section{Constrained fusion systems} \label{constrained}

We now look at a class of saturated fusion systems which have very simple, 
regular behavior:  the \emph{constrained} fusion systems. The main 
results here say that constrained fusion systems are always realized as 
fusion systems of finite groups in a predictable way, and have unique 
associated centric linking systems.

Let $\calf$ be an arbitrary saturated fusion system over a $p$-group $S$.  
Recall (Definition \ref{D:normal}) that a subgroup $Q\nsg{}S$ is 
\emph{normal in $\calf$} if each $\alpha\in\homf(P,P')$ extends to a 
morphism $\widebar{\alpha}\in\homf(PQ,P'Q)$ which sends $Q$ to itself.  
If $Q$ and $Q'$ are both normal in $\calf$, then clearly $QQ'$ is normal in 
$\calf$.  Hence, there is a unique maximal normal $p$-subgroup in $\calf$, 
which we denote $O_p(\calf)$ by analogy with the subgroup $O_p(G)$ of a 
finite group $G$. By Proposition~\ref{Q<|F}, $O_p(\calf)$ is contained in 
the intersection of all $\calf$-radical subgroups of $S$.  We are 
interested in the case when  $O_p(\calf)$ is itself $\calf$-centric, or 
equivalently, when there is a subgroup $P\nsg S$ which is both normal and 
centric in $\calf$.

\begin{Defi} \label{D:constrained}
A saturated fusion system $\calf$ over a $p$-group $S$ is 
\emph{constrained} if there is some $Q\nsg{}S$ which is $\calf$-centric 
and normal in $\calf$.
\end{Defi}

When $G$ is a finite $p'$-reduced group, then $G$ is said to be 
\emph{$p$-constrained} if there exists some normal $p$-subgroup $P\nsg G$ 
which is centric in $G$ (i.e., $C_G(P)\le{}P$).  (More generally, an 
arbitrary finite group $G$ is $p$-constrained if its $p'$-reduction 
$G/O_{p'}(G)$ is $p$-constrained.)  Our aim is to show that any 
constrained fusion system is the fusion system of a unique $p'$-reduced 
$p$-constrained group $G$.  This will be done by first showing that each 
constrained fusion system has a unique associated centric linking system 
$\call$, and then choosing $G$ to be a certain automorphism group in 
$\call$.

We first show that for any constrained fusion system, the obstruction 
groups to the existence and uniqueness of an associated centric linking 
system vanish. For any saturated fusion system $\calf$, let $\calz_\calf$ 
denote the functor on $\orb(\calf^c)$ defined by setting 
$\calz_\calf(P)=Z(P)$ for all $\calf$-centric $P\le{}S$.  (See 
\cite[\S3]{BLO2} for details.)

\begin{Prop}  \label{limz(constrained)}
Let $\calf$ be any constrained saturated fusion system over a $p$-group
$S$.  Then
    $$ \higherlim{\orb(\calf^c)}i(\calz_\calf) = 0
    \qquad\textup{for all $i>0$}. $$
In particular, there is a centric linking system $\call$ associated 
to $\calf$ which is unique up to isomorphism.
\end{Prop}

\begin{proof} Fix $Q\nsg{}S$ which is $\calf$-centric and normal in 
$\calf$. Let $P_1,P_2,\dots,P_m$ be $\calf$-conjugacy class 
representatives for all $\calf$-centric subgroups $P\le{}S$ such 
that$P\ngeq{}Q$, arranged such that $|P_i|\le|P_{i+1}|$ for each $i$.  For 
$i=0,1,\dots,m$, let $\calz_i\subseteq\calz_\calf$ be the subfunctor
	$$ \calz_i(P)= \begin{cases}  Z(P) & \textup{if $P$ is 
	$\calf$-conjugate to $P_j$ for some $j>i$} \\
	0 & \textup{otherwise.}
	\end{cases} $$
This gives a sequence of subfunctors 
$\calz_\calf\supseteq\calz_0\supseteq\calz_1\supseteq\cdots\supseteq 
\calz_m=0$, where for each $i=1,\dots,m$, $\calz_{i-1}/\calz_i$ vanishes 
except on subgroups $\calf$-conjugate to $P_i$.  Hence by 
\cite[Proposition~3.2]{BLO2}, 
	$$ \higherlim{\orb(\calf^c)}*(\calz_{i-1}/\calz_i) \cong
	\Lambda^*(\outf(P_i);Z(P_i)). $$
Furthermore, since $P_i\ngeq{}Q$, $N_{P_iQ}(P_i)/P_i\cong\Out_Q(P_i)$ is a 
nontrivial normal $p$-subgroup of $\outf(P_i)$ (normal by the same 
argument as the one used in the proof of Proposition \ref{Q<|F}), 
$\Lambda^*(\outf(P_i);Z(P_i))=0$ by \cite[Proposition~6.1(ii)]{JMO}.  This 
proves that $\higherlim{}*(\calz_i)=0$ for all $i$, and in particular that 
$\higherlim{}*(\calz_0)=0$.  Thus 
	\beq \higherlim{\orb(\calf^c)}*(\calz_\calf) \cong
	\higherlim{\orb(\calf^c)}*(\calz_\calf/\calz_0), \tag{1} \eeq
where $\calz_\calf/\calz_0$ is the quotient functor
	\beq (\calz_\calf/\calz_0)(P)= \begin{cases}  
	Z(P) = Z(Q)^P & \textup{if $P\ge Q$} \\
	0 & \textup{if $P\ngeq Q$.}
	\end{cases} \tag{2} \eeq

Now set $\Gamma=\outf(Q)$ and $S_0=\Out_S(Q)\cong{}S/Q$.  Thus 
$S_0\in\sylp{\Gamma}$.  Set $M=Z(Q)$, regarded as a 
$\zploc[\Gamma]$-module. Let $H^0M$ be the fixed-point functor on 
$\orb_{S_0}(\Gamma)$ defined by $H^0M(P)=M^P$. Then $H^0M$ is acyclic by 
\cite[Proposition~5.14]{JM} (shown more explicitly in 
\cite[Proposition~5.2]{JMO}).  So by (1), we will be done upon showing that 
	\beq \higherlim{\orb(\calf^c)}*(\calz_\calf/\calz_0) \cong
	\higherlim{\orb_{S_0}(\Gamma)}*(H^0M). \tag{3} \eeq
Since $Q$ is normal and centric in $\calf$, it is easy to check  that 
$\orb_{S_0}(\Gamma)$ is isomorphic to the full subcategory of 
$\orb(\calf^c)$ with objects the subgroups of $S$ containing $Q$. Under 
this identification, $H^0M$ is the restriction of $\calz_\calf/\calz_0$ by 
(2).  Isomorphism (3) now follows since $(\calz_\calf/\calz_0)(P)=0$ for 
all $P\ngeq{}Q$, and since there are no morphisms in $\orb(\calf^c)$ from 
an object in the subcategory to an object not in it.

The existence and uniqueness of a centric linking system associated to 
$\calf$ now follow from \cite[Proposition 3.1]{BLO2}.
\end{proof}

We are now ready to show that each constrained fusion system is the fusion 
system of a group. The following proposition includes Proposition 
\ref{D:constrained-intro}.

\begin{Prop}  \label{constrained1} 
Let $\calf$ be a constrained saturated fusion system over a $p$-group $S$. 
Then there is a unique finite $p'$-reduced $p$-constrained group $G$, 
containing $S$ as a Sylow $p$-subgroup, such that $\calf=\calf_S(G)$ as 
fusion systems over $S$.  Furthermore, if $\call$ is a centric linking 
system associated to $\calf$, then 
\begin{enumerate}  
\item $G\cong\Aut_\call(Q)$ for any subgroup
$Q\nsg{}S$ which is $\calf$-centric 
and normal in $\calf$; and 
\item $\call\cong\call_S^c(G)$.
\end{enumerate}
\end{Prop}

\begin{proof}  Using Proposition \ref{limz(constrained)}, fix a centric 
linking system $\call$ associated to $\calf$.  Let $\pi\:\call\rTo\calf^c$ 
denote the canonical projection functor.  By Lemma \ref{deltaPQ}, any 
choice of ``inclusion'' morphisms $\iota_P\in\Mor_\call(P,S)$ determines 
unique injections
	$$ \delta_{P,P'}\: N_S(P,P') \Right5{} \Mor_\call(P,P'), $$
for all $\calf$-centric subgroups $P,P'\le{}S$, which
satisfy the following conditions:
\begin{enumerate}\renewcommand{\labelenumi}{\textup{(\roman{enumi})}}%
\item $\pi(\delta_{P,P'}(g))=c_g\in\homf(P,P')$ for $g\in{}N_S(P,P')$; 
\item $\delta_{P,P}(g)=\delta_P(g)\in\Aut_\call(P)$ for $g\in{}P$; 
\item $\delta_{P,P''}(hg)=\delta_{P',P''}(h)\circ\delta_{P,P'}(g)$ for 
$g\in{}N_S(P,P')$ and $h\in{}N_S(P',P'')$; and
\item $\delta_{P,S}(1)=\iota_P$.
\end{enumerate}
Set $\iota_P^{P'}=\delta_{P,P'}(1)\in\Hom_\call(P,P')$ for all $P\le{}P'$ 
containing $Q$.  We think of these as the ``inclusion morphisms'' in 
$\call$.  By construction, $\iota_P^S=\iota_P$ and $\iota_P^P=\Id_P$ for 
all $P$, and $\iota_P^{P''}=\iota_{P'}^{P''}\circ\iota_P^{P'}$ whenever 
$P\le{}P'\le{}P''$.  

The proposition follows from the following points, which will be proven in 
Steps 1--2.
\begin{enumerate}\renewcommand{\labelenumi}{\textup{(\arabic{enumi})}}%
\item Assume $Q\nsg{}S$ is $\calf$-centric and normal in $\calf$, and 
$G=\Aut_\call(Q)$.  Then $G$ is $p'$-reduced and $p$-constrained; and we can 
identify $S$ with a subgroup of $G$ in such a way that $S\in\sylp{G}$ and 
$\calf=\calf_S(G)$.

\item Assume $G$ is $p'$-reduced and $p$-constrained, and such that 
$S\in\sylp{G}$ and $\calf=\calf_S(G)$.  Then $\call\cong\call_S^c(G)$.  
Also, if $Q\nsg{}S$ is any subgroup which is $\calf$-centric and normal in 
$\calf$, then $Q\nsg{}G$, and $G\cong\Aut_\call(Q)$.
\end{enumerate}

\smallskip

\noindent\textbf{Step 1: } Fix $Q\nsg{}S$ which is $\calf$-centric and 
normal in $\calf$, and set $G=\Aut_\call(Q)$.  Via the injection
	$$ \delta_{Q,Q}\:S=N_S(Q)\Right5{}\Aut_\call(Q)=G, $$
we identify $S$ as a subgroup of $G$.  Since $Q$ is fully normalized, 
	$$ S/Z(Q)\cong\Aut_S(Q) \in \sylp{\autf(Q)}, $$
where $\autf(Q)\cong{}G/Z(Q)$; and thus $S\in\sylp{G}$.  

Let $P,P'\le{}S$ be any pair of subgroups which contain $Q$.  For any 
$f\in\Mor_\call(P,P')$, there is (by Lemma \ref{L^n-prop}) a unique 
``restriction'' of $f$ to $Q$:  a unique element 
$\gamma(f)\in{}G=\Aut_\call(Q)$ such that 
$\iota_Q^{S}\circ{}\gamma(f)=f\circ\iota_Q^P$.  These restrictions clearly 
satisfy the following two conditions:
\begin{enumerate}\renewcommand{\labelenumi}{\textup{(\roman{enumi})}}%
\setcounter{enumi}{4}
\item $\gamma(f'\circ{}f)=\gamma(f'){\cdot}\gamma(f)$ for any 
$f'\in\Mor_\call(P',P'')$, any $Q\le{}P''\le{}S$; and
\item $\gamma(\delta_{P,P'}(x))=x$ for all $x\in{}N_S(P,P')$. 
\end{enumerate}
Furthermore, by condition (C) in Definition \ref{L-cat}, for each $g\in{}P$,
	$$ \delta_{S}(\pi(f)(g))\circ f=f\circ\delta_P(g) 
	\in\Mor_\call(P,S). $$
Upon restriction to $Q$ (and applying (v) and (vi)), this gives the relation
	$$ \delta_{Q,Q}(\pi(f)(g))\circ\gamma(f)=
	\gamma(f)\circ\delta_{Q,Q}(g) \in\Aut_\call(Q)=G. $$
In other words, under the identification 
$S=\delta_{Q,Q}(S)\le\Aut_\call(Q)=G$, this shows that
\begin{enumerate}\renewcommand{\labelenumi}{\textup{(\roman{enumi})}}%
\setcounter{enumi}{6}
\item $\gamma(f)\in{}N_G(P,P')$ and $c_{\gamma(f)}=\pi(f)\in\homf(P,P')$.
\end{enumerate}

Now,
	$$ C_G(Q) = \Ker\bigl[\Aut_\call(Q) \Right3{\pi} \autf(Q)\bigr] 
	= Z(Q) \,: $$
the first equality by (vii) (applied with $P=P'=Q$, so $\gamma(f)=f$), and 
the second by condition (A) in Definition \ref{L-cat}.  Thus $Q$ is 
centric in $G$. This also shows that $O_{p'}(G)=1$ (since 
$[O_{p'}(G),Q]=1$), and hence that $G$ is $p'$-reduced and $p$-constrained.

We must show that $\calf=\calf_S(G)$.  We first show that 
$\homf(P,P')\subseteq\Hom_G(P,P')$ for each $P,P'\le{}S$.  Since $Q$ is normal 
in $\calf$, each morphism in $\homf(P,P')$ extends to a morphism in 
$\homf(PQ,P'Q)$, and hence it suffices to work with subgroups $P,P'\ge{}Q$.  
In particular, $P$ and $P'$ are $\calf$-centric in this case.  For any 
$\varphi\in\homf(P,S)$, and any 
$f\in\Mor_\call(P,S)$ such that $\pi(f)=\varphi$, $\gamma(f)\in{}N_G(P,P')$ 
and $\varphi=c_{\gamma(f)}\in\Hom_G(P,P')$ by (vii), and thus
$\homf(P,P')\subseteq\Hom_G(P,P')$.  

Conversely, for any $P,P'\le{}S$ and any $g\in{}N_G(P,P')=N_G(PQ,P'Q)$, 
we claim that $c_g\in\homf(P,P')$.  Again, we can assume that
$P,P'\ge{}Q$.  Now, $c_g|_Q\in\autf(Q)$ by (vii) (applied with $P=P'=Q$ and 
$f=g$).  Since $Q=gQg^{-1}$ is $\calf$-centric, it is fully centralized in 
$\calf$, and so $c_g|_Q$ extends to an $\calf$-morphism defined on
	\[ N_{c_g|_Q} \defeq \{ x\in S\;|\; c_{gxg^{-1}} \in \Aut_S(Q)\}
	\ge P, \]
by condition (II) of Definition \ref{sat.Frob.}. In particular, $c_g|_Q$ 
extends to a morphism $\varphi\in\homf(P,S)\subseteq\Hom_G(P,S)$ (where the 
inclusion holds by the previous paragraph).  Let $h\in{}N_G(P,S)$ be such 
that $\varphi=c_h$.  Then $c_h|_Q=\varphi|_Q=c_g|_Q$, so $h=gx$ for some 
$x\in{}C_G(Q)$, and $C_G(Q)=Z(Q)$ as already shown.  Since $x\in P$, 
$c_x\in\autf(P)$, so $c_g\in\homf(P,S)$, and $c_g\in\homf(P,P')$ since 
$c_g(P)=gPg^{-1}\le{}P'$.

\smallskip

\noindent\textbf{Step 2: } Let $G$ be any finite $p'$-reduced 
$p$-constrained group such that $S\in\sylp{G}$ and $\calf=\calf_S(G)$.  
Then $\call\cong\call_S^c(G)$ by the 
uniqueness in Proposition \ref{limz(constrained)}.  

Let $Q\nsg{}S$ be any subgroup normal in $\calf=\calf_S(G)$.  Set 
$Q'=O_p(G)$; thus $C_G(Q')=Z(Q')$ by assumption.  Since $Q$ is normal in 
$\calf_S(G)$, for any $g\in{}G$, $c_g\in\Aut_G(Q')$ extends to some 
$c_{g'}\in\Aut_G(QQ')$; then $g^{-1}g'\in{}C_G(Q')=Z(Q')$, 
$g'\in{}N_G(QQ')$, and so $g\in{}N_G(QQ')$.  This shows that $QQ'\nsg{}G$, 
a normal $p$-subgroup, and hence $Q\le{}Q'=O_p(G)$.  Hence for any 
$g\in{}G$, $c_g\in\Aut_G(Q')$ restricts to an automorphism of $Q$ (since 
$Q$ is normal in $\calf_S(G)$), so $g\in{}N_G(Q)$, and this shows that 
$Q\nsg{}G$.  

In particular, if $Q$ is both $\calf$-centric and normal in $\calf$, then
	\beq \Aut_\call(Q) \cong \Aut_{\call_S^c(G)}(Q) \cong 
	N_G(Q)/O^p(C_G(Q))=G/1\cong G. \tag*{\proofbox} \eeq
\def\proofbox{}\end{proof}

It is in general not true, for a constrained fusion system $\calf$ over a 
$p$-group $S$ and a finite group $G$ such that $S\in\sylp{G}$ and 
$\calf=\calf_S(G)$, that $p$-subgroups of $S$ normal in $\calf$ are also 
normal in $G$.  For example, if $G=A_5$, $p=2$, $S\in\Syl_2(G)$, and 
$\calf= \calf_S(G)$, then $\calf$ \emph{is} a constrained fusion system, 
with $O_2(\calf)=S\cong{}C_2^2$.  Thus $S$ is normal in $\calf$, but not 
in $G$, in this case.  This shows the importance of assuming $G$ is 
$p'$-reduced and $p$-constrained.  In the given example, the unique 
$2'$-reduced $2$-constrained group associated to $\calf$ is $A_4$.  



\affiliationone{
C. Broto \& N. Castellana\\
Departament de Matem\`atiques\\ 
Universitat Aut\`onoma de Barcelona\\ 
E--08193 Bellaterra, Spain
\email{broto@mat.uab.es} \email{natalia@mat.uab.es}}
\affiliationtwo{
J. Grodal\\
Department of Mathematics\\ 
University of Chicago\\
Chicago, IL 60637, USA
\email{jg@math.uchicago.edu}}
\affiliationthree{
R. Levi\\
Department of Mathematical Sciences\\ 
University of Aberdeen, Meston Building 339\\
Aberdeen AB24 3UE, U.K.
\email{ran@maths.abdn.ac.uk}}
\affiliationfour{
B. Oliver\\
LAGA, Institut Galil\'ee\\
Av. J-B Cl\'ement\\ 
93430 Villetaneuse, France
\email{bob@math.univ-paris13.fr}}

\end{document}